\documentclass[leqno,a4paper,12pt]{book}
\usepackage{amssymb,amsmath,amsthm}
\usepackage[all]{xy}
\usepackage{color}
\usepackage{makeidx}

\numberwithin{equation}{section}

\usepackage{pgf,tikz}
\usetikzlibrary{arrows,calc,shapes,decorations.pathreplacing,patterns}
\usepackage{pgfplots}

\usepackage{enumerate}
\usepackage{mathrsfs}

\newtheorem{theorem}{Theorem}[section]
\newtheorem{proposition}[theorem]{Proposition}

\newtheorem{corollary}[theorem]{Corollary}

\theoremstyle{definition}

\newtheorem{definition}[theorem]{Definition}

\newtheorem{example}[theorem]{Example}

\newtheorem{remark}[theorem]{Remark}

\newcommand{\Rem}{\begin{remark}}
\newcommand{\enrem}{\end{remark}}
\newcommand{\on}{\operatorname}
\newcommand{\RR}{\mathrm{R}}

\newcommand{\C}{{\mathbb{C}}}

\newcommand{\R}{{\mathbb{R}}}

\newcommand{\Z}{{\mathbb{Z}}}

\def\phi{{\varphi}}
\def\epsilon{\varepsilon}

\newcommand{\cor}{{\bf k}}
\def\sha{\mathscr{A}}
\def\shb{\mathscr{B}}
\def\shc{\mathscr{C}}
\def\shd{\mathscr{D}}
\def\she{\mathscr{E}}
\def\shf{\mathscr{F}}

\def\shi{\mathscr{I}}

\def\shm{\mathscr{M}}
\def\shn{\mathscr{N}}
\def\sho{\mathscr{O}}

\def\shs{\mathcal{S}}

\def\shu{\mathscr{U}}

\newcommand{\rmpt}{{\rm pt}}

\newcommand{\into}{\hookrightarrow}

\renewcommand{\to}[1][]{\xrightarrow[]{#1}}
\newcommand{\from}[1][]{\xleftarrow[]{#1}}
\newcommand{\isoto}[1][]{\xrightarrow[#1]%
{{\raisebox{-.6ex}[0ex][-.6ex]{$\mspace{1mu}\sim\mspace{2mu}$}}}}

\newcommand{\To}{\mathop{\makebox[2em]{\rightarrowfill}}}

\newcommand{\muhom}{\mu hom}
\newcommand{\muHom}[1][]{\mathrm{Hom}^\mu_{\raise1.5ex\hbox to.1em{}#1}}
\newcommand{\Hom}[1][]{\mathrm{Hom}_{\raise1.5ex\hbox to.1em{}#1}}
\newcommand{\RHom}[1][]{\RR\mathrm{Hom}_{\raise1.5ex\hbox to.1em{}#1}}
\newcommand{\Ext}[2][]{\mathrm{Ext}_{\raise1.5ex\hbox to.1em{}#1}^{#2}}
\renewcommand{\hom}[1][]{{\mathscr{H}\mspace{-4mu}om}_{\raise1.5ex\hbox to.1em{}#1}}
\newcommand{\rhom}[1][]{{\RR\mathscr{H}\mspace{-3mu}om}_{\raise1.5ex\hbox to.1em{}#1}}
\newcommand{\rhomc}[1][]
{{\mathscr{H}\mspace{-3mu}om}^*_{\raise1.5ex\hbox to.1em{}#1}}

\newcommand{\ext}[2][]{{\mathscr{E}xt}_{\raise1.5ex\hbox to.1em{}#1}^{#2}}
\newcommand{\Tor}[2][]{\mathrm{Tor}^{\raise1.5ex\hbox to.1em{}#1}_{#2}}
\newcommand{\tens}[1][]{\mathbin{\otimes_{\raise1.5ex\hbox to-.1em{}{#1}}}}
\newcommand{\ltens}[1][]{\mathbin{\overset{\mathrm{L}}\tens}_{#1}}

\newcommand{\lltens}[1][]{{\mathop{\tens}\limits^{\rm L}}_{#1}}
\newcommand{\detens}{\underline{\etens}}

\newcommand{\etens}[1][]{{\mathbin{\boxtimes}}_{#1}}
\newcommand{\letens}[1][]{{\overset{\mathrm{L}}{\etens}_{#1}}}

\newcommand{\Endo}[1][]{\mathrm{End}_{\raise1.5ex\hbox to.1em{}#1}}

\newcommand{\Aut}[1][]{\mathrm{Aut}_{\raise1.5ex\hbox to.1em{}#1}}
\newcommand{\sect}{\Gamma}
\newcommand{\rsect}{\mathrm{R}\Gamma}

\newcommand{\conv}[1][]{\mathop{\circ}\limits_{#1}}

\newcommand{\aconv}[1][]{\mathop{\circ}\limits^{a}\limits_{#1}}

\newcommand{\sconv}[1][]{\mathop{\ast}\limits_{#1}}

\newcommand{\atimes}[1][]{\mathop{\times}\limits^{a}\limits_{#1}}

\newcommand{\oim}[1]{{#1}_*}
\newcommand{\eim}[1]{{#1}_!}

\newcommand{\roim}[1]{\RR{#1}_*}

\newcommand{\reim}[1]{\RR{#1}_!}
\newcommand{\opb}[1]{#1^{-1}}

\newcommand{\epb}[1]{#1^{!}}
\newcommand{\spb}[1]{#1^{*}}

\newcommand{\omDA}[1][M]{\omega_{\Delta_{#1}}}
\newcommand{\omDAI}[1][M]{\omega_{\Delta_{#1}}^{{\otimes-1}}}

\newcommand{\omAI}[1][M]{\omega_{{#1}}^{{\otimes-1}}}

\newcommand{\rc[1]}{{\Brc}({#1})}
\newcommand{\rcc[1]}{{\Brc}^c({#1})}

\newcommand{\md[1]}{{\rm Mod}({#1})}
\newcommand{\mdf[1]}{{\rm Mod^f}({#1})}

\newcommand{\mdcp[1]}{{\rm Mod^{{\rm c}}}({#1})}

\newcommand{\mdrc[1]}{{\rm Mod_{{\Rc}}}({#1})}

\newcommand{\gr}{{\rm gr}}
\newcommand{\Xsa}{{X_{\rm sa}}}
\newcommand{\Msa}{{M_{\rm sa}}}
\newcommand{\Msal}{{M_{\rm sal}}}

\newcommand{\rhosa}{\rho_{\rm sa}}
\newcommand{\rhosal}{\rho_{\rm sal}}

\newcommand{\tp}{{\rm tp}}

\newcommand{\Db}{{\mathcal D} b}
\newcommand{\Dbt}{{\mathcal D} b^\tp}
\newcommand{\thom}{thom}
\newcommand{\wtens}{\mathbin{\overset{\mathrm{w}}\tens}}
\newcommand{\Ot}[1][X]{\mathcal{O}^{\mathrm tp}_{#1}}
\newcommand{\Ow}[1][X]{\mathcal{O}^{\mathrm w}_{#1}}
\newcommand{\Oww}[1][X]{\mathcal{O}^{\omega}_{#1}}
\newcommand{\Ov}[1][X]{\Omega_{#1}}

\newcommand{\Cinft}{{\mathcal C}^{\infty,\tp}}
\newcommand{\Cin}[1][]{{\mathcal C}^{\infty,{#1}}}
\newcommand{\Cinf}{{\mathcal C}^{\infty}}

\newcommand{\sHH}{\mathscr{H}\mspace{-8mu}\mathscr{H}}
\newcommand{\HHO}[1][X]{\sHH({\sho}_{#1})}

\newcommand{\RHHOo}[2][]{\mathbb{HH}^0_{#1}(\sho_{#2})}

\newcommand{\HHE}[2][]{\mathscr{HH}_{#1}(\she_{#2})}

\newcommand{\RHHEo}[2][]{\mathbb{HH}^0_{#1}(\she_{#2})}

\newcommand{\sMH}{\mathscr{M}\mspace{-4mu}\mathscr{H}}
\newcommand{\MH}[2][]{\sMH_{#1}(\cor_{#2})}
\newcommand{\RMH}[2][]{\mathbb{MH}_{#1}({\cor}_{#2})}
\newcommand{\MHo}[2][]{\mathbb{MH}_{#1}^0({\cor}_{#2})}

\newcommand{\MHCo}[2][]{\mathbb{MH}_{#1}^0({\C}_{#2})}
\newcommand{\MHC}[2][]{\sMH_{#1}({\C}_{#2})}

\newcommand{\eqdot}{\mathbin{:=}}

\newcommand{\cl}{\colon}
\newcommand{\scbul}{{\,\raise.4ex\hbox{$\scriptscriptstyle\bullet$}\,}}

\newcommand{\tw}[1]{\widetilde{#1}}

\newcommand{\ol}{\overline}
\newcommand{\bl}{\bigl(}
\newcommand{\br}{\bigr)}
\newcommand{\lp}{{\rm(}}
\newcommand{\rp}{{\rm)}}

\newcommand{\Cc}{{\C\text{-c}}}
\newcommand{\Rc}{{\R\text{-c}}}
\newcommand{\cc}{{\text{cc}}}

\newcommand{\Sol}{{\shs\mspace{-2.5mu}\mathit{ol}}}
\newcommand{\Solt}{{\shs\mspace{-2.5mu}\mathit{ol}^t}}
\newcommand{\Solo}{{\shs\mspace{-2.5mu}\mathit{ol}^0}}
\newcommand{\Solot}{{\shs\mspace{-2.5mu}\mathit{ol}^{0,t}}}
\newcommand{\GSol}{{\mathrm{Sol}}}

\newcommand{\spec}{{\mathrm{spec}}}
\newcommand{\ba}{\begin{array}}
\newcommand{\ea}{\end{array}}

\newcommand{\bnum}{\begin{enumerate}[{\rm(i)}]}
\newcommand{\enum}{\end{enumerate}}
\newcommand{\banum}{\begin{enumerate}[{\rm(a)}]}
\newcommand{\eanum}{\end{enumerate}}

\newcommand{\eq}{\begin{eqnarray}}
\newcommand{\eneq}{\end{eqnarray}}
\newcommand{\eqn}{\begin{eqnarray*}}
\newcommand{\eneqn}{\end{eqnarray*}}

\newcommand{\Proof}{\begin{proof}}
\newcommand{\QED}{\end{proof}}
\newcommand{\Prop}{\begin{proposition}}
\newcommand{\enprop}{\end{proposition}}

\def\rop{{\rm op}}

\def\Op{{\rm Op}}
\def\dist{{\rm dist}}

\def\hh{{\rm hh}}
\def\mueu{{\mu\rm eu}}
\def\hol{{\rm hol}}
\newcommand{\DQ}{\ensuremath{\mathrm{DQ}}}

\newcommand{\HK}{\ensuremath{\mathrm{TK}}}
\newcommand{\Ch}{\ensuremath{\mathrm{Ch}}}
\newcommand{\Td}{\ensuremath{\mathrm{Td}}}

\DeclareMathOperator{\id}{id}

\DeclareMathOperator{\supp}{supp}
\DeclareMathOperator{\ori}{or}
\DeclareMathOperator{\chv}{char}

\newcommand{\Supp}{\on{Supp}}
\newcommand{\Der}[1][]{\mathsf{D}^{#1}}
\newcommand{\Derb}{\Der[\mathrm{b}]}

\newcommand{\SSi}{\mathrm{SS}}
\newcommand{\SSid}{\mathrm{SS}_\Delta}

\newcommand{\RD}{\mathrm{D}}
\newcommand{\RDO}{\mathrm{D}_\sho}
\newcommand{\RDD}{\mathrm{D}_\shd}

\newcommand{\Int}{{\rm Int}}
\newcommand{\coh}{{\rm coh}}

\newcommand{\dT}{{\dot{T}}}

\newcommand{\indlim}[1][]{\mathop{\varinjlim}\limits_{#1}}
\newcommand{\sindlim}[1][]{\smash{\mathop{\varinjlim}\limits_{#1}}\,}

\newcommand{\sprolim}[1][]{\smash{\mathop{\varprojlim}\limits_{#1}}\,}

\newcommand{\inddlim}[1][]{\mathop{``{\varinjlim}"}\limits_{#1}}
\newcommand{\sinddlim}[1][]{\smash{\mathop{``{\varinjlim}"}\limits_{#1}}\,}

\newcommand{\hs}{\hspace*}

\newcommand{\bwr}{\mbox{\large{$\wr$}}}

\newcommand{\Cd}{\mathrm{C}}
\newcommand{\tK}{\widetilde{K}}

\newcommand{\htd}{{\rm td}}
\newcommand{\HOD}[1][X]{\mathcal{HD}({\sho}_{#1})}

\newcommand{\II[1]}{{\rm I}({#1})}
\newcommand{\IIrc[1]}{{\rm I_{\R-c}}({#1})}

\newcommand{\indc}{{\rm Ind}(\shc)}

\newcommand{\roimdpim}{{R\dot\pi_{M*}}}

\newcommand{\shbD}{\shb_\Delta}
\newcommand{\shcD}{\shc_\Delta}
\newcommand{\shbDU}{\shb_\Delta^\vee}
\newcommand{\shcDU}{\shc_\Delta^\vee}

\def\sal{{\rm sal}}

\def\rhosal{\rho_{\sal}}

\def\Msal{{M_{\sal}}}

\begin{document}
\author{Pierre Schapira}
\title{Three lectures on\\Algebraic Microlocal Analysis.\\
{\small Spring school,
Northwestern University, May 2012}}
\date{}
\maketitle
\tableofcontents

%\begin{tikzpicture}
%\draw (1,0)--(0,1)--(-1,0)--(0,-1) --cycle;
%\end{tikzpicture}
%\begin{tikzpicture}
%\draw (0,0) circle (1cm);
%\end{tikzpicture}
%\begin{tikzpicture}
%\draw (1,0)--(0,1)--(-1,0)--(0,-1) --cycle;
%\end{tikzpicture}

\chapter[Microlocalization]{Microlocalization of sheaves}
\noindent
{\bf Abstract.}
This first talk is a survey talk with some historical comments and I refer to~\cite{Sc10} for a more detailed overview.

I will first explain the notions of Sato's hyperfunctions and microfunctions, at the origin of the story, 
and I will describe the Sato's  microlocalization functor which was first motivated by problems of Analysis (see~\cite{SKK73}). Then I will briefly recall the main features of the microlocal theory of sheaves of~\cite{KS90} with
 emphasize on the functor $\muhom$ which will be the main tool for the second talk.

${}$
  
\section{Generalized functions}
In the sixties, people were used to work with various spaces of generalized functions constructed with the tools of functional analysis. Sato's construction of hyperfunctions in 59-60 is at the opposite of this practice: he uses 
purely algebraic tools and complex analysis. The importance of Sato's definition is twofold:
first, it is purely algebraic (starting with the analytic object $\sho_X$), and second it highlights  the link between real and complex geometry. 
(See \cite{Sa60} and see~\cite{Sc07} for an exposition of Sato's work.)

Consider first the case where
$M$ is an open subset of the real line $\R$ and let $X$ be an open neighborhood 
of $M$ in the complex line $\C$ satisfying $X\cap\R=M$.
The space $\shb(M)$ of hyperfunctions on $M$ is given by
\eqn
&&\shb(M)=\sho(X\setminus M)/\sho(X).
\eneqn
It is easily proved, using the solution of the Cousin problem,  that this
space depends only on $M$, not on the choice of $X$, and that the 
correspondence 
$U\mapsto \shb(U)$ ($U$ open in $M$)  defines a  flabby sheaf $\shb_M$ on $M$. 

With Sato's definition, the boundary values always exist and are no more a
limit in any classical sense. 
\begin{example}
(i) The Dirac function at $0$ is 
\eqn
&&\delta(0)=\frac{1}{2i\pi}(\frac{1}{x-i0}-\frac{1}{x+i0}).
\eneqn
Indeed, if $\phi$ is a $C^0$-function on $\R$ with compact support, one has
\eqn
&&\phi(0)=\lim_{\epsilon\to[>]0}\frac{1}{2i\pi}\int_\R(\frac{\phi(x)}{x-i\epsilon}-\frac{\phi(x)}{x+i\epsilon})dx.
\eneqn
(ii) The holomorphic function $\exp(1/z)$ defined on $\C\setminus\{0\}$ has a boundary value as a 
hyperfunction (supported by $\{0\}$) not as a distribution. 
\end{example}

On a  real analytic manifold $M$ of dimension $n$, the sheaf $\shb_M$ was originally defined as
\eqn
&&\shb_M=H^n_M(\sho_X)\tens \ori_{M}
\eneqn
where $X$ is a complexification of $M$ and $\ori_M$ is the orientation sheaf on $M$. 
It is shown that this object  is concentrated in degree $0$. 
Since $X$ is oriented, 
Poincar{\'e}'s duality gives the isomorphism $\RD'_X(\C_{M})\simeq \ori_{M}\,[-n]$ 
(see~\eqref{eq:dual1} below for the definition of $\RD_M'$). An equivalent definition of 
hyperfunctions is thus given  by
\eq\label{eq:hyp}
\shb_M=\rhom[\C_X](\RD'_X(\C_{M}),\sho_X).
\eneq
Let us define the notion of  ``boundary value'' in this settings.
Consider a subanalytic open subset $\Omega$ of $X$ and denote by
$\overline\Omega$ its closure. Assume that:
\eqn
&&\left\{\begin{array}{l}
\RD'_X(\C_{\Omega})\simeq \C_{\overline\Omega},\\
M\subset \overline\Omega.
\end{array}\right.
\eneqn
The morphism $\C_{\overline\Omega}\to\C_M$ defines by duality the morphism
$\RD'_X(\C_{M})\to\RD'_X(\C_{\overline\Omega})\simeq \C_{\Omega}$. Applying the functor
$\RHom(\scbul,\sho_X)$, we get 
the boundary value morphism
\eq\label{eq:bvmor}
&&{\rm b}\cl \sho(\Omega)\to \shb(M).
\eneq
When considering operations on hyperfunctions such as integral transforms, one is naturally lead to consider more general sheaves of generalized functions such as 
$\rhom(G,\sho_X)$ where $G$ is a constructible sheaf. We shall come back on this point later.

Similarly as in dimension one, we can represent the sheaf $\shb_M$ by using \v{C}ech cohomology of coverings of $X\setminus M$. For example, let $X$ be a Stein open subset of $\C^n$ and set $M=\R^n\cap X$.
Denote by $x$ the coordinates on $\R^n$ and by $x+iy$ the coordinates on $\C^n$.
One can recover $\C^n\setminus\R^n$ by $n+1$ open half-spaces
$V_i=\langle y,\xi_i\rangle>0$ ($i=1,\dots,n+1$). For $J\subset\{1,\dots,n+1\}$ set $V_J=\bigcap_{j\in J}V_j$. 
Assuming $n>1$, we have the isomorphism $H^n_M(X;\sho_X)\simeq H^{n-1}(X\setminus M;\sho_X)$. 
Therefore, setting $U_J=V_J\cap X$, one has
\eqn
&&\shb(M)\simeq \sum_{\vert J\vert=n}\sho_X(U_J)/\sum_{\vert K\vert=n-1}\sho_X(U_K).
\eneqn
On a real analytic manifold $M$, any hyperfunction $u\in\sect(M;\shb)$ is a (non unique) sum of boundary values of holomorphic functions defined in tubes with edge $M$. Such a decomposition 
leads to the so-called Edge of the Wedge theorem and was intensively studied in the seventies (see~\cite{Mr67, BI73}).

Then comes naturally the following problem: how to recognize the directions associated with these tubes?
The answer is given by the Sato's microlocalization functor.

\section{Microlocalization}

Unless otherwise specified, all manifolds are real, say of class $C^\infty$ and
$\cor$ denotes  a commutative  unital ring with finite global homological dimension.

We denote by $\cor_M$ the constant sheaf on $M$ with stalk $\cor$,
by $\Derb(\cor_M)$ the bounded derived category of  sheaves of $\cor$-modules on $M$ and by 
$\Derb_{\cc}(\cor_{M})$ the full triangulated subcategory of 
$\Derb(\cor_{M})$ consisting of  cohomologically constructible objects. If $M$ is real analytic, we denote by $\Derb_{\Rc}(\cor_{M})$
the triangulated category of $\R$-constructible sheaves.

We denote by $\omega_M$ the dualizing complex on $M$. Then
$\omega_M\simeq\ori_M\,[\dim M]$ where $\ori_M$ is the orientation sheaf and $\dim M$ the dimension of 
$M$. We shall use the duality functors 
\eq\label{eq:dual1}
&&\RD_M'F=\rhom(F,\cor_M),\quad \RD_M F=\rhom(F,\omega_M).
\eneq

For a locally closed subset $A$ of $M$, we denote by $\cor_{MA}$ the sheaf which is the constant sheaf on $A$ with stalk $\cor$ and which is $0$ on $M\setminus A$. If there is no risk of confusion, we simply denote it by $\cor_A$.

\subsubsection*{Fourier-Sato transform}
The classical Fourier transform interchanges (generalized) functions on a vector space $V$ and 
(generalized) functions on the dual vector space $V^*$. The idea of extending this formalism 
to sheaves, hence to replacing an isomorphism of spaces with an equivalence of categories,  seems to have 
appeared first in  Mikio Sato's construction of microfunctions in the 70's.

Let $\tau:E\to M$ be a finite dimensional real vector bundle over a real manifold $M$ with fiber
dimension $n$ and let $\pi: E^*\to M$ be the dual vector bundle. 
Denote by $p_1$ and $p_2$ the first and second projection
defined on $E\times_M E^*$, and define:
\eqn
 P=\{(x,y)\in E\times_M E^*; \langle x,y\rangle \geq 0\},\\
 P'=\{(x,y)\in E\times_M E^*; \langle x,y\rangle \leq 0\}.
\eneqn
Consider the diagram:
\eqn
&&\xymatrix{
&E\times_M E^* \ar[ld]_-{p_1}\ar[rd]^-{p_2} &\\
E\ar[rd]_-\tau&&E^*\ar[ld]^-\pi\\
&M.&
}\eneqn
Denote by $\Derb_{\R^+}(\cor_E)$ the full triangulated subcategory  of $\Derb(\cor_E)$
consisting of conic sheaves, that is, objects with locally constant cohomology on the orbits of $\R^+$.

\begin{definition}
Let $F\in\Derb_{\R^+}(\cor_E)$, $G\in \Derb_{\R^+}(\cor_{E^*})$. One sets:
\eqn
F^\wedge &=& \reim{p_2}(\opb{p_1}F)_{P'}\simeq  \roim{p_2}(\rsect_P\opb{p_1}F),\\
G^\vee & =&   \roim{p_1}(\rsect_{P'}\epb{p_2}G)\simeq\reim{p_1}(\epb{p_2}G)_{P}.
\eneqn
\end{definition}
The main result of the theory is the following.
\begin{theorem}\label{th:Fourier}
The two functors $(\cdot)^\wedge$ and
$(\cdot)^\vee$ are inverse to each other, hence define an equivalence of categories
$\Derb_{\R^+}(\cor_E) \simeq \Derb_{\R^+}(\cor_{E^*})$.
\end{theorem}
\begin{example}
(i) Let $\gamma$ be a closed proper convex cone in $E$ with $M\subset \gamma$. Then:
$$(\cor_\gamma)^\wedge \simeq \cor_{\Int \gamma^\circ}.$$
Here $\gamma^\circ$ is the polar cone to $\gamma$, a
closed convex cone in $E^*$ and $\Int \gamma^\circ$ denotes its interior.

\noindent
(ii) Let $\gamma$ be an open convex cone in $E$. Then: 
$$(\cor_\gamma)^\wedge \simeq \cor_{\gamma^{\circ a}}\tens \ori_{E^*/M}\, [-n].$$
Here $\lambda^a = -\lambda$, the image of $\lambda$ by the antipodal map. 

\noindent
(iii) Let $(x)=(x',x'')$ be coordinates on $\R^n$ with $(x')=(x_1,\dots,x_p)$ and $(x'')=(x_{p+1},\dots,x_n)$.
Denote by $(y)=(y',y'')$ the dual coordinates on $(\R^n)^*$. 
Set 
\eqn
&&\gamma=\{x;x'^2-x''^2\geq0\},\quad \lambda=\{y;y'^2-y''^2\leq0\}.
\eneqn
Then $(\cor_\gamma)^\wedge\simeq \cor_\lambda[-p]$. (See~\cite{KS97}.)
\end{example}

\subsubsection*{Specialization}
Let $\iota\cl N\hookrightarrow M$ be the  embedding of a closed submanifold $N$ of $M$. Denote by 
$\tau_M\cl T_NM\to N$ the normal bundle to $N$. 

If $F$ is a sheaf on $M$, its restriction to $N$, denoted $F\vert_N$, may be viewed as a global object, namely the 
direct image by $\tau_M$ of a sheaf $\nu_MF$ on $T_NM$, called the specialization of $F$ along $N$. 
Intuitively, $T_NM$ is the set of light rays issued from $N$ in $M$ and 
the germ of $\nu_NF$ at a normal vector $(x;v)\in T_NM$ is the germ at $x$ of the restriction of $F$ along the light ray $v$.

One constructs a new manifold $\tw M_N$, called the normal deformation of $M$ along $N$,  together with the maps
\eq\label{diag:specializ}
\xymatrix{
T_NM\ar[r]^-s\ar[d]_{\tau_M}&\tw M_N\ar[d]^-p&\Omega\ar[l]_-j\ar[ld]^-{\tw p}\\
N\ar[r]_-{\iota}&M
}, \quad t\cl \tw M_N\to\R, \,\,\Omega=\{\opb{t}(\R_{>0})\}
\eneq
with the following properties.
Locally, after choosing a local coordinate system $(x',x'')$ on $M$ such that $N=\{x'=0\}$, we have 
$\tw M_N=M\times\R$, $t\cl \tw M_N\to\R$ is the projection, $\Omega=\{(x;t)\in M\times\R;t>0\}$,
$p(x',x'',t)=(tx',x'')$, $T_NM=\{t=0\}$.

Let $S\subset M$ be a locally closed subset. The Whitney normal cone $C_N(S)$ is a closed conic subset of $T_NM$ given by
\eqn
&&C_N(S)=\ol{\opb{\tw p}(S)}\cap T_NM
\eneqn
where, for a set $A$,  $\ol A$ denotes the closure of $A$.
One defines the specialization functor 
\eqn
&&\nu_N\cl\Derb(\cor_M)\to\Derb(\cor_{T_NM})
\eneqn
by a similar formula, namely:
\eqn
&&\nu_NF\eqdot \opb{s}\oim{j}\opb{\tw p}F.
\eneqn
Clearly, $\nu_NF\in\Derb_{\R^+}(\cor_{T_NM})$, that is, $\nu_NF$  is a conic sheaf for the $\R^+$-action on $T_NM$. Moreover,
\eqn
&&\roim{\tau_M}\nu_NF\simeq \nu_NF\vert_N\simeq F\vert_N.
\eneqn
For an open cone $V\subset T_NM$, one has 
\eqn
&&H^j(V;\nu_NF)\simeq\sindlim[U]H^j(U;F)
\eneqn
where $U$ ranges through the family of open subsets of $M$ such that 
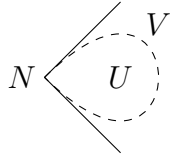
\eqn
&&C_N(M\setminus U)\cap V=\emptyset. 
\eneqn
\begin{center}
\begin{tikzpicture}
\draw (0,0) -- (1,1);
\draw (0,0) -- (1,-1);
\draw [dashed] (0,0) ..controls +(2,2) and  +(2,-2)  .. (0,0);
%\fill[color=gray](0,0) ..controls +(1.9,1.9) and  +(1.9,-1.9)  .. (0,0);
%\draw [->] (0,0) -- (-1,0);
\draw (-0.3,0) node {$N$};
\draw (1.5,0.7) node {$V$};
\draw (1.0,0) node {$U$};
\end{tikzpicture}
\end{center}
\subsubsection*{Microlocalization}
Denote by 
$\pi_M\cl T^*_NM\to N$ the conormal bundle to $N$ in $M$, that is, the dual bundle to $\tau_M\cl T_NM\to N$.

If $F$ is a sheaf on $M$,  the sheaf of sections of $F$ supported by $N$, denoted $\rsect_NF$, may be viewed as a global object, namely the 
direct image by $\pi_M$ of a sheaf $\mu_MF$ on $T^*_NM$. 
Intuitively, $T^*_NM$ is the set of ``walls'' (half-spaces)  in $M$ passing through $N$ and 
the germ of $\mu_NF$ at a conormal vector $(x;\xi)\in T^*_NM$ is the germ at $x$ of the sheaf of sections of $F$ supported by closed tubes with edge $N$ and which are almost 
the half-space associated with $\xi$.

More precisely, the microlocalization of $F$ along $N$, denoted $\mu_NF$, is the Fourier-Sato transform of 
$\nu_NF$, hence is an object of $\Derb_{\R^+}(\cor_{T^*_NM})$. It satisfies:
\eqn
&&\roim{\pi_M}\mu_NF\simeq \mu_NF\vert_N\simeq \rsect_NF.
\eneqn
For a convex open cone $V\subset T^*_NM$, one has 
\eqn
&&H^j(V;\mu_NF)\simeq \indlim[U,Z]H^j_{U\cap Z}(U;F),
\eneqn
where $U$ ranges over the family of open subsets of $M$ such that $U\cap N=\pi_M(V)$ and 
$Z$ ranges over the family of closed subsets of $M$ such that $C_M(Z)\subset V^\circ$ where $V^\circ$ is the polar cone to $V$. 
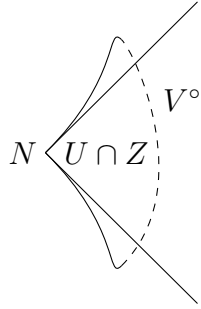
\begin{center}
\begin{tikzpicture}
\draw (0,0) ..controls +(1,1) and  +(-0.2,0.2)  .. (1,1.5);
\draw (0,0) ..controls +(1,-1) and  +(-0.2,-0.2)  .. (1,-1.5);
\draw [dashed] (1,-1.5) ..controls +(1,1) and  +(0.2,-0.2)  .. (1,1.5);
%\fill[color=gray](1,-1.5) ..controls +(0.9,0.9) and  +(0.19,-0.19)  .. (1,1.5);
\draw (0,0) -- (2,2);
\draw (0,0) -- (2,-2);
%\draw [->] (0,0) -- (-1,0);
\draw (-.3,0) node {$N$};
\draw (1.8,0.7) node {$V^\circ$};
\draw (0.8,0) node {$U\cap Z$};
\end{tikzpicture}
\end{center}

\subsubsection*{Back to hyperfunctions}
Assume now that $M$ is a real analytic manifold and $X$ is a complexification of $M$. 
First notice the isomorphisms 
\eqn
&&M\times_XT^*X\simeq \C\tens_\R T^*M\simeq T^*M\oplus\sqrt{-1}T^*M.
\eneqn
One deduces the  isomorphism 
\eq\label{eq:sqrt1}
&&T^*_MX\simeq \sqrt{-1}T^*M. 
\eneq
The sheaf $\shc_M$ on $T^*_MX$ of Sato's microfunction (see~\cite{SKK73}) is defined as
\eqn
&&\shc_M\eqdot\mu_M(\sho_X)\tens\opb{\pi_M}\omega_M.
\eneqn
It is shown that this object  is concentrated in degree $0$. 
Therefore, we have an isomorphism
\eqn
&&\spec\cl \shb_M\isoto\oim{\pi_M}\shc_M
\eneqn
and Sato defines the analytic wave front set of a hyperfunction $u\in\sect(M;\shb_M)$ as the support of 
$\spec(u)\in\sect(T^*_MX;\shc_M)$. 

Consider a closed convex proper cone $Z\subset T^*_MX$ which contains the zero-section $M$.
Then  $\spec(u)\subset Z$ if and only if $u$ is the boundary value of a holomorphic function 
 defined in a tuboid $U$ with profile the interior of the polar tube to $Z^a$ 
 (where $Z^a$ is the image of $Z$ by the antipodal map), that is, satisfying
 \eqn
 &&C_M(X\setminus U)\cap \Int Z^{\circ a}=\emptyset.
 \eneqn
Moreover, the sheaf $\shc_M$ is conically flabby. Therefore, any hyperfunction may be decomposed as a sum of boundary values of holomorphic functions $f_i$'s defined in suitable tuboids $U_i$ and 
if we have hyperfunctions 
$u_i$ ($i=1,\dots N$) satisfying $\sum_ju_j=0$, there exist hyperfunctions $u_{ij}$ ($i,j=1,\dots N$) 
such that 
\eqn
&&u_{ij}=-u_{ji}, \quad u_i=\sum_ju_{ij} \mbox{ and }\spec(u_{ij})\subset\spec(u_i)\cap\spec(u_j).
\eneqn
When translating this result in terms of boundary values of holomorphic functions, we get the so-called ``Edge of the wedge theorem'', already mentioned. 

Sato's introduction of the sheaf $\shc_M$ was the starting point of an intense activity in the domain of linear partial differential equations after H\"ormander adapted Sato's ideas to  classical analysis with the help of the (usual) Fourier transform. See~\cite{Ho83} and also~\cite{BI73, Sj82} for related constructions. 
Note that the appearance of $\sqrt{-1}$ in the usual Fourier transform may be understood 
as following from the isomorphism \eqref{eq:sqrt1}.

\section{Microsupport}
The microsupport of sheaves has been introduced in~\cite{KS82} and developed in~\cite{KS85,KS90}.
Roughly speaking, the microsupport of $F$ describes the codirections of non propagation of $F$. 
The idea of microsupport takes its origin in the study of linear PDE and particularly in the study of hyperbolic systems.

\begin{definition}
Let $F\in \Derb(\cor_M)$ and let $p\in T^*M$. 
One says that $p\notin\SSi(F)$ if there exists an open neighborhood
$U$ of $p$ such that for any $x_0\in M$ and any
real $\Cd^1$-function $\phi$ on $M$ defined in a neighborhood of $x_0$ 
with $(x_0;d\phi(x_0))\in U$, one has
$(\rsect_{\{x;\phi(x)\geq\phi(x_0)\}} F)_{x_0}\simeq0$.
\end{definition}
In other words, $p\notin\SSi(F)$ if the sheaf $F$ has no cohomology 
supported by ``half-spaces'' whose conormals are contained in a 
neighborhood of $p$. 
\begin{itemize}
\item
By its construction, the microsupport is $\R^+$-conic, that is,
invariant by the action of  $\R^+$ on $T^*M$. 
\item
$\SSi(F)\cap T^*_MM=\pi_M(\SSi(F))=\Supp(F)$.
\item
The microsupport satisfies the triangular inequality:
if $F_1\to F_2\to F_3\to[{\;+1\;}]$ is a
distinguished triangle in  $\Derb(\cor_M)$, then 
$\SSi(F_i)\subset\SSi(F_j)\cup\SSi(F_k)$ for all $i,j,k\in\{1,2,3\}$
with $j\not=k$. 
\end{itemize}

\begin{example}\label{ex:microsupp}
(i) If $F$ is a non-zero local system on $M$ and $M$ is connected, then $\SSi(F)=T^*_MM$.

\noindent
(ii) If $N$ is a closed submanifold of $M$ and $F=\cor_N$, then 
$\SSi(F)=T^*_NM$, the conormal bundle to $N$ in $M$.

\noindent
(iii) Let $\phi$ be a $\Cd^1$-function such that $d\phi(x)\not=0$ 
whenever $\phi(x)=0$.
Let $U=\{x\in M;\phi(x)>0\}$ and let $Z=\{x\in M;\phi(x)\geq0\}$. 
Then 
\eqn
&&\SSi(\cor_U)=U\times_MT^*_MM\cup\{(x;\lambda d\phi(x));\phi(x)=0,\lambda\leq0\},\\
&&\SSi(\cor_Z)=Z\times_MT^*_MM\cup\{(x;\lambda d\phi(x));\phi(x)=0,\lambda\geq0\}.
\eneqn
\end{example}
For a precise definition of being co-isotropic (one also says involutive),
we refer to~\cite[Def.~6.5.1]{KS90}.

\begin{theorem}\label{th:coisotr}
Let $F\in \Derb(\cor_M)$. Then its microsupport 
$\SSi(F)$ is co-isotropic.
\end{theorem} 

Assume now that $(X,\sho_X)$ is a complex manifold and denote as usual by $\shd_X$ the sheaf of rings of finite order differential operators on $X$. 
For a coherent $\shd_X$-module $\shm$, one denotes by $\chv(\shm)$ its characteristic variety, a closed conic 
complex analytic subvariety of $T^*X$. One also sets for short
\eqn
&&\Sol(\shm)\eqdot\rhom[\shd](\shm,\sho_X).
\eneqn
After identifying $X$  with its real underlying manifold,
the link between the microsupport of sheaves and the characteristic
variety of coherent $\shd$-modules is given  by
\begin{theorem}\label{th:ssinchar1}
Let $\shm$ be a coherent $\shd$-module. Then 
$\SSi(\Sol(\shm))=\chv(\shm)$.
\end{theorem}
The inclusion $\SSi(\Sol(\shm))\subset \chv(\shm)$
is the most useful in practice. Its proof  only makes use of
the Cauchy-Kowalevsky theorem in its precise form given by Petrovsky and 
Leray (see~\cite[\S~9.4]{Ho83})
and of purely algebraic arguments.
As a corollary of Theorems~\ref{th:coisotr} and~\ref{th:ssinchar1}, one recovers the fact that the characteristic variety 
of a coherent $\shd_X$-module is co-isotropic, a theorem of~\cite{SKK73} which also have a purely algebraic proof due to Gabber~\cite{Ga81}.

\section{The functor $\muhom$}

We denote by $\delta\cl M\to M\times M$ the diagonal embedding and we set $\Delta=\delta(M)$. For short, we also denote by $\delta$ the isomorphism
\eqn
&&\delta\cl T^*M\isoto T^*_\Delta(M\times M),\quad (x;\xi)\mapsto (x,x;\xi,-\xi).
\eneqn

Let us briefly recall the main properties of the functor $\muhom$, a variant of Sato's microlocalization functor. 
\eqn
&&\muhom\cl \Derb(\cor_M)^\rop\times \Derb(\cor_M)\to\Derb(\cor_{T^*M}),\\
&&\muhom(G,F)\eqdot\opb{\delta}\mu_\Delta\rhom(\opb{q_2}G,\epb{q_1}F)
\eneqn
where $q_i$ ($i=1,2$) denotes the $i$-th projection on  $M\times M$. Note that
\eqn
&&\roim{\pi_M}\muhom(G,F)\simeq\rhom(G,F),\\
&&\muhom(\cor_N,F)\simeq\mu_N(F)\mbox{ for $N$ a closed submanifold of $M$},\\
&&\supp\muhom(G,F)\subset\SSi(G)\cap\SSi(F),\\
&&\muhom(G,F)\simeq\mu_\Delta(F\letens \RD_M G)\mbox{ if $G$ is constructible.}
\eneqn
In some sense, $\muhom$ is the sheaf of microlocal morphisms. More precisely, for $p\in T^*M$, we have;
\eqn
&&H^0\muhom(G,F)_p\simeq \Hom[\Derb(\cor_M;p)](G,F)
\eneqn
where the category  $\Derb(\cor_M;p)$ is  the localization of $\Derb(\cor_M)$ by the subcategory of sheaves whose microsupport does not contain $p$. 

There is an interesting phenomena which holds with $\muhom$ and not with $\rhom$.  Indeed, assume $M$ is real analytic. Then, although the category $\Derb_{\Rc}(\cor_{M})$ of constructible sheaves 
does not admit a Serre functor, it admits 
a kind of microlocal Serre functor, as shown by the isomorphism, functorial with respect to $F$ and $G$ 
(see~\cite[Prop.~8.4.14]{KS90}):
\eqn
&&\RD_{T^*M}\muhom(F,G)\simeq\muhom(G,F)\tens\opb{\pi_M}\omega_M.
\eneqn
This confirms the fact that to fully understand constructible sheaves, it is natural to look at them 
microlocally, that is, in $T^*M$. This is also in accordance with the ``philosophy'' of Mirror Symmetry which interchanges the category of coherent $\sho_X$-modules on a complex manifold $X$ with the Fukaya category on a symplectic manifold $Y$. In case of $Y=T^*M$, the Fukaya category is equivalent to the category of $\R$-constructible sheaves on $M$, according to Nadler-Zaslow~\cite{Na09,NZ09} (see also~\cite{FLTZ10} for related results.)

\section{An application: elliptic pairs}
Denote by $\dT^*M$ the set $T^*M\setminus T^*_MM$ and denote by 
$\dot\pi_M$ the restriction of $\pi_M\cl T^*M\to M$ to $\dT^*M$. 
If $H\in\Derb_{\R^+}(\cor_{T^*M})$ is   a conic sheaf on $T^*M$, there is the Sato's distinguished triangle
\eqn
&&\reim{\pi_M}H\to\roim{\pi}H\to\roimdpim H\to[+1].
\eneqn
Applying this result with $H=\muhom(G,F)$ and assuming $G$ is constructible, we get the 
distinguished triangle
\eqn
&&\RD_M'G\tens F\to\rhom(G,F)\to\roimdpim\muhom(G,F).
\eneqn

\begin{theorem}{\rm (The Petrovsky theorem for sheaves.)}\label{th:petrowshv}
Assume that $G$ is constructible and $\SSi(G)\cap\SSi(F)\subset T^*_MM$. Then the natural morphism
\eqn
&&\rhom(G,\cor_M)\tens F\to\rhom(G,F)
\eneqn
is an isomorphism.
\end{theorem}

Let us apply this result when $X$ is a complex manifold.
For $G\in\Derb_\Rc(\C_X)$, set
\eqn
&&\sha_G=\sho_X\tens G, \quad \shb_G\eqdot\rhom(\RD_X'G,\sho_X).
\eneqn
Note that if $X$ is the complexification of a real analytic manifold $M$ and we choose $G=\C_M$, we 
recover the sheaf of real analytic functions and the sheaf of hyperunctions:
\eqn
&&\sha_{\C_M}=\sha_M,\quad \shb_{\C_M}=\shb_M.
\eneqn
Now let $\shm\in\Derb_\coh(\shd_X)$. 
According to~\cite{ScSn94}, one says that the pair $(G,\shm)$ is elliptic if 
$\chv(\shm)\cap\SSi(G)\subset T^*_XX$. 

\begin{corollary}{\rm \cite{ScSn94}}   \label{cor:petrowshv}
Let $(\shm,G)$ be an elliptic pair. 
\banum
\item
We have the canonical isomorphism:
\eq\label{eq:solasolb}
&&\rhom[\shd_X](\shm,\sha_G)\isoto\rhom[\shd_X](\shm,\shb_G).
\eneq
\item
Assume moreover that $\Supp(\shm)\cap\Supp(G)$ is compact.
Then the cohomology of the complex 
$\RHom[\shd_X](\shm,\sha_G)$
is finite dimensional.
\eanum
\end{corollary}
To prove the part~(b) of the corollary, one represents  the left hand side of
the global sections of \eqref{eq:solasolb}
by a complex of topological vector spaces of type DFN and the right hand side 
by a complex of topological vector spaces of type FN.

  \chapter[Microlocal Euler classes]{Microlocal Euler classes and Hochschild homology}
{\bf Abstract.}
This is a joint work with Masaki Kashiwara announced in~\cite{KS12}.
On a complex manifold $(X,\sho_X)$,  the 
Hochschild homology is a powerful tool to construct characteristic classes of coherent modules and to get index theorems. 
Here, I will show how to adapt this formalism to a wide class of sheaves on a real manifold $M$ by using the functor $\muhom$ of microlocalization.
This construction applies in particular to 
constructible sheaves on real manifolds and $\shd$-modules on complex manifolds, or more generally to elliptic pairs.

\section[Complex manifolds]{Hochschild homology on complex manifolds}
Hochschild homology of $\sho$-modules has given rise to a vast literature. Let us quote 
in particular~\cite{Hu06,Ca05,CaW07,Ra08}.

Consider a complex manifold $(X,\sho_X)$ and denote by
$\omega_X^\hol$ the dualizing complex in the category of
$\sho_X$-modules, that is, $\omega_X^\hol=\Omega_X\,[d_X]$, where
$d_X$ is the complex dimension of $X$ and $\Omega_X$ is the sheaf of 
holomorphic forms of degree $d_X$. 
We shall use the classical six operations for $\sho$-modules, $\spb{f}$, $\roim{f}$, $\epb{f}$, $\reim{f}$,
$\ltens[\sho]$ and $\rhom[\sho]$. 
In particular we have the two duality functors
\eqn
&&\RDO'(\scbul)=\rhom[\sho_X](\scbul,\sho_X),\\
&&\RDO(\scbul)=\rhom[\sho_X](\scbul,\omega^\hol_X)
\eneqn
as well as the external product that we denote by $\letens[\sho]$.
Denote by $\delta\cl X \hookrightarrow X\times X$ the diagonal embedding
and let $\Delta=\delta(X)$. We set
\eq
\sho_\Delta\eqdot\oim{\delta}\sho_X,\quad
\omega_X^{\hol,\tens-1}\eqdot\RDO'\omega_X^{\hol}, \quad
\omega_\Delta^{\hol,\tens-1}\eqdot\oim{\delta}\omega_X^{\hol,\tens-1}.
\eneq
It is well-known that 
\eq\label{eq:omegadd}
&&\omega_\Delta^{\hol,\tens-1}
\simeq\rhom[\sho_{X\times X}](\sho_\Delta,\sho_{X\times X}).
\eneq
The Hochschild homology of $\sho_X$ is usually defined by   
\eq\label{eq:hh1}
&&\HHO[X]=\opb{\delta}\bl\sho_\Delta\ltens[\sho_{X\times X}]\sho_\Delta\br§.
\eneq
Note the isomorphisms
\eqn
&&\xymatrix{
&\HHO[X]\ar[ld]^-\sim\ar[rd]_-\sim&\\
\spb{\delta}\oim{\delta}\sho_X\ar[rr]^-\sim&&\epb{\delta}\eim{\delta}\omega^\hol_X
}
\eneqn
and the canonical  isomorphisms
\eqn
\spb{\delta}\oim{\delta}\sho_X&\simeq&\opb{\delta}\rhom[\sho_{X\times X}]\bl\omega_\Delta^{\hol,\tens-1},\sho_\Delta\br,\\
\epb{\delta}\eim{\delta}\omega^\hol_X&\simeq& \opb{\delta}\rhom[\sho_{X\times X}]\bl\sho_\Delta,\omega_\Delta^{\hol}\br.
\eneqn
For a closed  subset $S$ of $X$, we set:
\eq\label{eq:hh5}
&&\RHHOo[S]{X}=H^0(X;\rsect_S\HHO[X]).
\eneq

Let $\shf\in\Derb_\coh(\sho_X)$. The morphisms $\RDO'\shf\ltens[\sho_X]\shf\to\sho_X$ and $\RDO\shf\ltens[\sho_X]\shf\to\omega^\hol_X$
give by adjunction the morphisms 
\eqn
&&\RDO'\shf\letens[\sho_{X\times X}]\shf\to \sho_\Delta ,\quad  \RDO\shf\letens[\sho_{X\times X}]\shf \to\omega^\hol_\Delta
\eneqn
and then by duality the morphisms
\eqn
&&\omega^{\hol,\tens-1}_\Delta\to \RDO'\shf\letens[\sho_{X\times X}]\shf \to \sho_\Delta,\quad   \sho_\Delta\to\RDO\shf\letens[\sho_{X\times X}]\shf\to 
\omega^\hol_\Delta
\eneqn
and the composition defines the Hochschild classes of $\shf$:
\eq\label{eq:HHO}
&&\hh_\sho(\shf)\in H^0_{\supp(\shf)}(X;\opb{\delta}\oim{\delta}\sho_X), \quad
\tw\hh_\sho(\shf)\in H^0_{\supp(\shf)}(X;\epb{\delta}\eim{\delta}\omega_X).
\eneq

One can compose Hochschild homology and the Hochschild class commutes with the composition of kernels. 
More precisely, consider complex manifolds  $X_i$ ($i=1,2,3$). 
\begin{itemize}
\item
We write $X_{ij}\eqdot X_i\times X_j$ ($1\leq i,j\leq3$),
$X_{123}=X_1\times X_2\times X_3$,
$X_{1223}=X_1\times X_2 \times X_2\times X_3$, etc.  
\item
We denote by $q_i$ the
projection $X_{ij}\to X_i$ or the projection $X_{123}\to X_i$ and by $q_{ij}$
the projection $X_{123}\to X_{ij}$. 
\end{itemize}
Let $K_{ij}\in\Derb_\coh(\sho_{X_{ij}})$  \lp$i=1,2$, $j=i+1$\rp. One sets
\eqn
&&K_{12}\conv[2]K_{23}=\reim{q_{13}}(\spb{q_{12}} K_{12}\ltens[\sho_{X_{123}}]\spb{q_{23}} K_{23})
\eneqn 

\begin{theorem}\label{th:functhhO}
\banum
\item
There is a natural morphism
\eqn
&&\HHO[X_{12}]\conv[2]\HHO[X_{23}]\to\HHO[X_{13}].
\eneqn
\item
Let $S_{ij}\subset X_{ij}$ be a closed subset \lp$i=1,2$, $j=i+1$\rp. Assume that 
$q_{13}$ is proper over $S_{12}\times_{X_2} S_{23}$ and set $S_{13}=q_{13}(S_{12}\times_{X_2}S_{23})$.
Then the morphism above induces a map
\eqn
&&\conv[2]\cl \RHHOo[S_{12}]{X_{12}}\tens \RHHOo[S_{23}]{X_{23}} \to\RHHOo[S_{13}]{X_{13}}.
\eneqn
\item
Let  $K_{ij}$ be as above and assume that $\supp(K_{ij})\subset S_{ij}$.
Set  $K_{13}=K_{12}\conv[2]K_{23}$ and $\tw K_{13}=(K_{12}\tens\omega^{\hol\tens-1}_{2})\conv[2]K_{23}$. Then 
$K_{13}$ and $\tw K_{13}$ belong to $\Derb_\coh(\sho_{X_{13}})$ and we have the equalities in 
$\RHHOo[S_{13}]{X_{13}}$:
 \eqn
\hh_\sho(K_{13})=\hh_\sho(K_{12})\conv[2] \hh_\sho(K_{23}),\quad
\tw\hh_\sho(\tw K_{13})=\tw\hh_\sho(K_{12})\conv[2] \tw\hh_\sho(K_{23}).
\eneqn
\eanum
\end{theorem}

This theorem shows in particular that 
 the Hochschild class commutes with external product, inverse image and proper direct image.

Theorem~\ref{th:functhhO} seems to be  well-known from the specialists although it is difficult to find a precise statement
(see however~\cite{CaW07,Ra10}).
The construction of the Hochschild homology as well as Theorem~\ref{th:functhhO} (including complete proofs)
have been extended when replacing $\sho_X$ with a
so-called $\DQ$-algebroid stack $\sha_X$ in~\cite{KS12b}.

Coming back to $\sho_X$-modules, the  Hodge cohomology of $\sho_X$ is given by:
\eq\label{def:HOD}
&&\HOD[X]\eqdot\bigoplus_{i= 0}^{d_X}\Omega_X^i\,[i],
\text{ an object of }\Derb(\sho_X).
\eneq
There is a commutative diagram constructed by Kashiwara in~\cite{Ka91}  in which $\alpha_X$ is the HKR (Hochschild-Kostant-Rosenberg) isomorphism and $\beta_X$ is a kind of dual HKR isomorphism:
\eq\label{diag:RR1}
&&\xymatrix{
\spb{\delta}\oim{\delta}\sho_X\ar[d]^-\sim_-{\alpha_X}\ar[rr]^-\sim_{\htd}
                                         &&\epb{\delta}\eim{\delta}\omega^\hol_X\\
\HOD[X]\ar[rr]^-\sim_{\tau}&&\HOD[X]\ar[u]^-\sim_-{\beta_X}.
}\eneq
If $\shf\in\Derb_\coh(\sho_X)$, the Chern character of $\shf$ is the image by $\alpha_X$ of $\hh_\sho(\shf)$.

In~\cite{Ka91}  Kashiwara made the conjecture that the arrow $\tau$ 
making the diagram commutative is given by the cup product by the Todd class of $X$.
This conjecture has recently been proved by Ramadoss~\cite{Ra08} in the algebraic case (after preliminary important results by Markarian)  and Grivaux~\cite{Gr09} in the analytic case (and with a very simple proof).
Since the morphism $\beta_X$ commutes with proper direct images, we get a new and functorial approach to the 
Riemann-Roch-Hirzebruch-Grothendieck theorem.

\section[Sheaves]{Microlocal homology}
We keep the notations of Lecture I. In particular $\omega_M$ denotes the dualizing complex on $M$ and
 $\RD'_M$ is the duality functor. We set
\eq
\omega_\Delta\eqdot\oim{\delta}\omega_M,\quad
\omega_M^{\tens-1}\eqdot \RD'_M\omega_M, \quad\omega_\Delta^{\tens-1}\eqdot\oim{\delta}\omega_M^{\tens-1}.
\eneq
Let $M_i$ ($i=1,2,3$) be manifolds. 
\begin{itemize}
\item
For short, we write as above 
$M_{ij}\eqdot M_i\times M_j$ ($1\leq i,j\leq3$),
$M_{123}=M_1\times M_2\times M_3$, etc.
\item
We will often write for short $\cor_i$ instead of $\cor_{{M_i}}$ and $\cor_{\Delta_i}$ 
instead of $\cor_{\Delta_{M_i}}$, $\pi_i$ instead of $\pi_{M_i}$, etc.
\item
We denote by $q_i$ the
projection $M_{ij}\to M_i$ or the projection $M_{123}\to M_i$ and by $q_{ij}$
the projection $M_{123}\to M_{ij}$. Similarly, we denote by $p_i$ the
projection $T^*M_{ij}\to T^*M_i$ or the projection $T^*M_{123}\to T^*M_i$ and
by $p_{ij}$ the projection $T^*M_{123}\to T^*M_{ij}$. 
\item 
We also need tointroduce the maps $p_{j^a}$ or $p_{ij^a}$, the composition of $p_{j}$ or $p_{ij}$ and the antipodal
map on $T^*M_j$. 
\end{itemize}

We consider the operations of composition of kernels. For $K_{ij}\in\Derb(\cor_{M_{ij}})$ ($i=1,2$, $j=i+1$), we set
\eqn\label{eq:conv}
&&\ba{rcl}
K_1\conv[2] K_2&\eqdot&\reim{q_{13}}\opb{\delta_2}(K_1\etens K_2)
      \simeq \reim{q_{13}}  (\opb{q_{12}} K_1\tens\opb{q_{23}}),\\
K_1\sconv[2] K_2&\eqdot&\roim{q_{13}}\bl\epb{\delta_2}(K_1\etens K_2)\tens \opb{q_{2}}\omega_{2}\br.
\ea
\eneqn
We have  a natural
morphism $K_1 \conv K_2  \to K_1 \sconv K_2$. It is an isomorphism if
$p_{12^a}^{-1}\SSi(K_1)\cap p_{23^a}^{-1}\SSi(K_2)\to T^*M_{13}$ is proper.

We also define the composition of kernels on cotangent bundles. For $L_{i}\in\Derb(\cor_{T^*_{M_{ij}}})$ 
($i=1,2$, $j=i+1$), we set

\eqn\label{eq:aconv}
\hs{-0ex}\ba{rcl}
L_1\aconv[2] L_2&\eqdot&\reim{p_{13^a}}(\opb{p_{12^a}} L_1\tens\opb{p_{23^a}} L_2).
\ea
\eneqn
For $K_1,F_1 \in \Derb(\cor_{M_{12}})$ and $K_2,F_2 \in \Derb(\cor_{M_{23}})$
there exists a canonical morphism:
\eq\label{eq:microcompk}
&&\muhom(K_1,F_1)\aconv[2]\muhom(K_2,F_2)\to\muhom(K_1\sconv[2]K_2,F_1\conv[2]F_2).
\eneq

We also define the corresponding operations for subsets of cotangent bundles.
Let $A\subset T^*M_{12}$ and $B\subset T^*M_{23}$. We set $A\aconv[2] B=p_{ 13}(A\atimes[2]B)$
where $A\atimes[2]B=\opb{ p_{12^a}}(A)\cap\opb{ p_{23}}(B)$.

If there is no risk of confusion, we simply denote by $\delta^a$ the map:
\eqn
&&
\delta^a\cl \xymatrix@C=2.5ex{T^*M\ar@{^{(}->}[r]& T^*(M\times M)},\quad (x;\xi)\mapsto (x,x;\xi,-\xi).
\eneqn

\begin{definition}\label{def:muHH} 
Let $\Lambda$ be a closed conic subset of $T^*M$. We set  
\eqn
\MH[]{M}&\eqdot&\opb{(\delta^a)}\muhom(\cor_{\Delta_M},\omDA[M]),\\
\MHo[\Lambda]{M}&\eqdot& H^0_\Lambda(T^*M;\MH[]{M}).
\eneqn
We call  $\MH[]{M}$ the microlocal homology of $M$.
\end{definition}
We have isomorphisms
\eqn
\MH[]{M}&\simeq&
\opb{(\delta^a)}\mu_\Delta(\omega_\Delta)\simeq\opb{\pi_M}\omega_M
\eneqn
and the isomorphism $\MH[]{M}\simeq\opb{\pi_M}\omega_M$
plays the role of the HKR isomorphism in the complex case.

We have the analogue of Theorem~\ref{th:functhhO}~(a) and~(b). (For the part (c), see Theorem~\ref{th:HH1} below.)

Let $i=1,2$, $j= i+1$ and let $\Lambda_{ij}$ be a closed conic subset of $T^*M_{ij}$.
Assume that
\eq\label{hyp:112233proper0}
&&\Lambda_{12}\atimes[2]\Lambda_{23} \mbox{ is proper over }T^*M_{13}.
\eneq
Note that this hypothesis is equivalent to 
\eqn&&
\left\{\parbox{63ex}{$
\opb{p_{12^a}}(\Lambda_{12})\cap\opb{p_{23^a}}(\Lambda_{23})\cap(T^*_{M_1}M_1\times T^*M_2\times T^*_{M_3}M_3)\subset T^*_{M_{123}}M_{123},\\
q_{13} \mbox{ is proper on }\pi_{12}(\Lambda_{12})\times_{M_2}\pi_{23}(\Lambda_{23}).
$}\right.
\eneqn
Set 
\eq\label{eq:lambda13}
&&\Lambda_{13}=\Lambda_{12}\aconv[2]\Lambda_{23}.
\eneq

\begin{theorem}\label{th:functhhShv}
\banum
\item
There is a natural morphism
\eq\label{eq:morhh1}
&&\MH[]{M_{12}} \conv[2]\MH[]{M_{23}}\to\MH[]{M_{13}}.
\eneq
\item
Let $\Lambda_{ij}\subset T^*M_{ij}$ be as above and assume~\eqref{hyp:112233proper0}.
Then morphism~\eqref{eq:morhh1}    induces a map
\eq\label{eq:morhh2}
&&\conv[2]\cl \MHo[\Lambda_{12}]{M_{12}}\tens \MHo[\Lambda_{23}]{X_{23}} \to\MHo[\Lambda_{13}]{X_{13}}.
\eneq
\eanum
\end{theorem}
The construction of the morphism~\eqref{eq:morhh1} uses~\eqref{eq:microcompk}, which makes 
the computations not easy. Fortunately, we have the following result.

\begin{proposition}\label{prop:complagcyc}
Let $M_i$ \lp$i=1,2,3$\rp\, be manifolds and let $\Lambda_{ij}$ be a closed conic subset of $T^*M_{ij}$ \lp $ij=12,13,23$\rp. We have a commutative diagram 
\eq\label{HKRreal}
&&\ba{l}\xymatrix{
\MH[]{12}\aconv[2]\MH[]{23}\ar[r]\ar[d]^{\bwr}& \MH[]{13}                       \ar[d]^{\bwr}\\
 \pi_{{12}}^{-1}\omega_{{12}}\aconv[2]\pi_{{23}}^{-1}\omega_{{23}}\ar[r]&\pi_{{13}}^{-1}\omega_{{13}}.
}\ea\label{dia:muconv}
\eneq
\end{proposition}

Here the bottom horizontal arrow is induced by
\eqn
&&\opb{p_{12^a}}\opb{\pi_{{12}}}\omega_{{12}}\tens \opb{p_{23^a}}\opb{\pi_{{23^a}}}\omega_{{23}}
\simeq \opb{\pi_{1}}\omega_{1}\etens 
\omega_{T^*M_2}\etens\opb{\pi_{3}}\omega_{3}\\
&&\text{  and}\\
&&\reim{p_{13^a}}\bl \pi_{1}^{-1}\omega_{1}\etens 
\omega_{T^*M_2}\etens\pi_{M_3}^{-1}\omega_{3}\br
\To \pi_{1}^{-1}\omega_{1}\etens \pi_{3}^{-1}\omega_{3}.
\eneqn

\begin{remark}
(i) If we consider that the isomorphism $\MH[]{M}\simeq\opb{\pi}\omega_M$
is a real analogue of the Hochschild-Kostant-Rosenberg isomorphism, then the commutativity of Diagram~\eqref{HKRreal} says that, contrarily to the complex case, the real HKR isomorphism commutes with inverse and direct images.

\vspace{0.3ex}\noindent
(ii) As a particular case of Proposition~\ref{prop:complagcyc}, we get canonical isomorphisms
\eqn
&&\MH{M}\tens\MH{M}\simeq \opb{\pi}\omega_M\tens \opb{\pi}\omega_M\simeq\omega_{T^*M}.
\eneqn
Hence, $\MH{M}$ behaves as a ``square root'' of the dualizing complex.
\end{remark}
  
${}$
  
\section[Microlocal Euler classes]{Trace kernels and microlocal Euler classes}
A trace kernel $(K,u,v)$ on $M$ is the data of $K\in\Derb(\cor_{M\times M})$ together with 
morphisms $(u,v)$
\eqn
&&\cor_\Delta\to[u] K\to[v]\omega_\Delta.
\eneqn
Setting $\SSid(K)\eqdot\SSi(K)\cap T^*_\Delta (M\times M)$,
the morphism $u$ gives an element 
of $H^0_{\SSid(K)}(T^*M;\muhom(\cor_\Delta,K))$ whose image by $v$ 
is  the microlocal Euler class of $K$ 
\eqn
&&  \mueu_M(K)\in \MHo[\SSid(K)]{M})\simeq H^0_{\SSid(K)}(T^*M;\opb{\pi}\omega_M).
\eneqn

If $M=\rmpt$, a Hochschild kernel $K$ is nothing but an object of $\Derb(\cor)$ together with linear maps $\cor\to K\to\cor$. The composition gives the element $\mueu(K)$ of $\cor$. If $\cor$ is a field of characteristic zero
 and $K=L\tens L^*$ where $L$ is a bounded complex of $\cor$-modules with finite dimensional cohomology and $L^*$ is its dual, one recovers
the classical Euler-Poincar\'e index of $L$, that is,  $\mueu(K)=\chi(L)$.

Let $i=1,2$, $j= i+1$ and let $\Lambda_{iijj}$ be a closed conic subset of $T^*M_{iijj}$.
Assume that
\eq\label{hyp:112233proper}
&&\Lambda_{1122}\atimes[22]\Lambda_{2233} \mbox{ is proper over }T^*M_{1133}.
\eneq
Set $\Lambda_{1133}=\Lambda_{1122}\aconv[22]\Lambda_{2233}$ and $\Lambda_{ij}=\Lambda_{iijj}\cap T^*_{\Delta_{ij}}M_{iijj}$.

\begin{theorem}\label{th:HH1}
Let $K_{ij}$ be a trace kernel on $M_{ij}$ with $\SSi(K_{ij})\subset\Lambda_{iijj}$.
Assume~\eqref{hyp:112233proper}, set $\tK_{23}=\omDAI[2]\conv[2]K_{23}
\simeq(\omAI[2]\letens\cor_{233})\ltens K_{23}$  and set
$K_{13}=K_{12}\conv[22]\tK_{23}$. Then 
\banum
\item
$K_{13}$ is a trace kernel on $M_{13}$,
\item
 $\mueu_{M_{13}}(K_{13})=\mueu_{M_{12}}(K_{12})\aconv[2] \mueu_{M_{23}}(K_{23})$
as elements of $\MHo[\Lambda_{13}]{13}$.
\eanum
\end{theorem}

As an application, one can perform the external product, the proper direct image and the non characteristic inverse image of trace kernels and compute their microlocal Euler classes. 

Consider in particular the case where $\Lambda_{1}$ and $\Lambda_{2}$ are two closed conic subsets of $T^*M$ satisfying
 the transversality condition
\eq\label{eq:transvlambda}
&&\Lambda_{1}\cap\Lambda_{2}^a\subset T^*_MM.
\eneq
Then applying Theorem~\ref{th:HH1} and composing the external product   with the restriction to the diagonal, 
we get a  convolution map:
\eqn
\star\cl\RMH[\Lambda_1]{M}\times\RMH[\Lambda_2]{M}\to\RMH[\Lambda_1+\Lambda_2]{M}.
\eneqn
\begin{proposition}\label{pro:convmueu}
 Let  $K_i$ be a trace kernels with $\SSi_\Delta(K_{i})\subset\Lambda_{i}$ \lp$i=1,2$\rp\,
 and assume~\eqref{hyp:112233proper}. Then the object $K_1\ltens (\cor_M\letens\omAI)\ltens K_2$ is a trace kernel on $M$ and 
\eqn
&&\mueu_M(K_1\ltens(\cor_M\letens\omAI)\ltens K_2)=\mueu_M(K_1)\star\mueu_M(K_2).
\eneqn
In particular if $\supp K_1\cap\supp K_2$ is compact, we have
\eqn
\mueu\bl\rsect(M\times M;K_1\ltens (\cor_M\letens\omAI)\ltens K_2)\br&=&\int_M(\mueu(K_1)\star\mueu(K_2))\vert_M\\
&=&\int_{T^*M}\mueu(K_1)\cup\mueu(K_2).
\eneqn
\end{proposition}
We shall apply this result to elliptic pairs. 

${}$

\section[Constructible sheaves]{Microlocal Euler class of constructible sheaves}
Let us denote by  $\Derb_\cc(\cor_M)$ the full  triangulated  subcategory of 
$\Derb(\cor_M)$ consisting of cohomologically constructible  sheaves
and let $G\in\Derb_\cc(\cor_M)$. 

The evaluation morphism 
$G\ltens\RD_M G\to\omega_M$ gives by adjunction the morphism $G\letens\RD_M G\to\omega_\Delta$.
By duality, one gets the morphism $\cor_\Delta\to G\letens\RD_M G$. To summarize, we have 
the morphisms in  $\Derb_\cc(\cor_{M\times M})$:
\eq
&&\cor_\Delta\to G\letens\RD G\to\omega_\Delta.
\eneq
Denote by $\HK(G)$ the Hochschild kernel so constructed. If $G$ is $\R$-constructible, the class 
$\mueu_M(\HK(G))$ is nothing but the Lagrangian cycle of $G$ constructed by Kashiwara~\cite{Ka85}.
In the sequel, if there is no risk of confusion,  we simply denote this class by $\mueu_M(G)$.

One
recovers the classical functorial properties of Lagrangian cycles. 
Let $f\cl M\to N$ be a morphism of manifolds.
To $f$ one associates the maps
\eqn
&&T^*M\from[f_d] M\times_{N}T^*N\to[f_\pi] T^*N
\eneqn
%\eq\label{eq:fdfpi}
%&&\ba{c}\xymatrix@C=10ex{
%T^*M\ar[dr]_-{\pi_{M}}
%                 &M\times_{N}T^*N\ar[d]^\pi\ar[l]_-{f_d}\ar[r]^-{f_\pi}
%                                        & T^*N\ar[d]^-{\pi_{N}}\\
%                  &M\ar[r]^-f   &N.
%}\ea\eneq
There are natural morphsim
\eqn
&&f_\mu\cl \eim{f_\pi}\opb{f_d}\opb{\pi_M}\omega_M\to \opb{\pi_N}\omega_N,\\
&&f^\mu\cl \eim{f_d}\opb{f_\pi}\opb{\pi_N}\omega_N\to \opb{\pi_M}\omega_M.
\eneqn

\begin{itemize}
\item
Let $F\in\Derb_\Rc(\cor_M)$ and assume $f$ is proper on $\supp(F)$, or equivalently, $f_\pi$ is proper on $\opb{f_d}\SSi(F)$. Then 
$\mueu(\roim{f}F)=f_\mu\mueu(F)$,
\item
Let $G\in\Derb_\Rc(\cor_N)$ and assume that $f$ is non characteristic for $G$, that is, $f_d$ is proper on 
$\opb{f_\pi}\SSi(G)$. Then $\mueu(\opb{f}G)=f^\mu\mueu(D)$.
\end{itemize}

\section[$\shd$-modules]{Microlocal Euler class of  $\shd$-modules}

In this section, we denote by $X$ a complex manifold of complex dimension $d_X$ and the base ring $\cor$ is the field $\C$. One denotes by $\shd_X$ the sheaf of $\C_X$-algebras of (finite order) holomorphic differential operators on $X$ and  refer to~\cite{Ka03} for a detailed exposition of the theory of $\shd$-modules. 

We also denote by $\Derb_\coh(\shd_X)$ the full triangulated subcategory of $\Derb(\shd_X)$ consisting of objects with coherent cohomology.
We denote by $\RDD\cl \Derb(\shd_X)\to\Derb(\shd_X)$ the duality functor for left $\shd$-modules:
\eqn
&&\RDD\shm\eqdot\rhom[\shd_X](\shm,\shd_X)\tens[\sho_X]\omega_X^{\hol,\tens-1}.
\eneqn
We denote by $\scbul\detens\scbul$ the external product for $\shd$-modules:
\eqn
&&\shm\detens\shn\eqdot \shd_{X\times X}\tens[\shd_X\etens\shd_X](\shm\etens\shn).
\eneqn
Let $\Delta$ be the diagonal of $X\times X$. 
The left $\shd_{X\times X}$-module $H^{d_X}_{[\Delta]}(\sho_{X\times X})$ (the algebraic cohomology with support in $\Delta$) is denoted as usual by $\shbD$. We also introduce
$\shbDU\eqdot \shbD\,[2d_X]$.
For  a coherent $\shd_X$-module $\shm$, we have the isomorphism
\eqn
\rhom[\shd_X](\shm,\shm)  &\simeq&\rhom[\shd_{X\times X}](\shbD,\shm\detens\RDD\shm)\,[d_X].
\eneqn
We get the morphisms
\eq\label{eq:hhforD}
&&\shbD\to\shm\detens\RDD\shm\,[d_X]\to\shbDU
\eneq
where the second morphism is deduced by duality. 

Denote by $\she_{T^*X}$ the sheaf on $T^*X$ of microdifferential operators of ~\cite{SKK73}. For a coherent $\shd_X$-module $\shm$ set
\eqn
&&\shm^E\eqdot \she_{T^*X}\tens[\opb{\pi}\shd_X]\opb{\pi}\shm.
\eneqn
Recall that, denoting by  $\chv(\shm)$ the characteristic variety of $\shm$, we have $\chv(\shm)=\supp(\shm^E)$. 
Set
\eqn
&&\shcD\eqdot \shbD^E,\quad \shcDU\eqdot (\shbDU)^E.
\eneqn
Let $\Lambda$ be a closed conic subset of $T^*X$.
One sets 
\eqn
&&\HHE[]{T^*X}=\opb{(\delta^a)}\rhom[\she_{X\times X}](\shcD,\shcDU),\\
&&\RHHEo[\Lambda]{T^*X}=H^0_\Lambda(T^*X;\HHE[]{T^*X}).
\eneqn
One calls $\HHE[]{T^*X}$ the Hochschild homology of $\she_{T^*X}$. 

We deduce from~\eqref{eq:hhforD} the morphisms
\eq  \label{eq:EHHmor0}
&&\shcD\to (\shm\detens\RDD\shm)^E\,[d_X]\to\shcDU  
\eneq
which define the Hochschild class of $\shm$:
\eq\label{eq:hhM}
&&\hh_\she(\shm)\in\RHHEo[\chv(\shm)]{T^*X}.
\eneq
We shall make a link between the Hochschild class of $\shm$ and the microlocal Euler class of a Hochschild kernel attached to the sheaf of holomorphic solutions of $\shm$. 
We have
\eqn
&&\Omega_{X\times X}\,[-d_X]\ltens[\shd_{X\times X}]\shbD\simeq\C_\Delta,\\
&&\Omega_{X\times X}\,[-d_X]\ltens[\shd_{X\times X}] \shbDU\simeq\omega_\Delta.
\eneqn
Now remark that for $\shn_1,\shn_2\in\Derb_\coh(\shd_X)$, we have a natural morphism
\eqn\label{eq:homDhomC}
&&\rhom[\opb{\pi}\shd_X](\opb{\pi}\shn_1,\shn_2^E)\to\muhom(\Omega_X\ltens[\shd_X]\shn_1,\Omega_X\ltens[\shd_X]\shn_2).
\eneqn
One deduces the morphisms
\eqn
\rhom[\she_{X\times X}](\shcD,\shcDU)
  &&\simeq\rhom[\opb{\pi}\shd_{X\times X}](\opb{\pi}\shbD,(\shbDU)^E)\\
 &&\hspace{-7ex}\to \muhom(\Omega_{X\times X}
 \ltens[\shd_{X\times X}]\shb_\Delta{}^{\tens-1},\Omega_{X\times X}\ltens[\shd_{X\times X}]\shb_\Delta)\\
 &&\hspace{-7ex}\simeq\muhom(\C_\Delta,\omega_\Delta).
\eneqn
Since  all the arrows above are isomorphisms, we get
\eqn
&&\HHE[]{T^*X}\simeq\MHC[]{X}.
\eneqn
Recall that the Hochschild homology of $\she_{T^*X}$ has been already calculated in~\cite{BG87}.

By this isomorphism, $\hh_\she(\shm)$ belongs to  $\MHCo[\chv(\shm)]{X}$ and this class coincides with that already introduced in~\cite{ScSn94}.  

Applying the functor $\Omega_{X\times X}[-d_X]\ltens[\shd_{X\times X}]\scbul$ to the morphisms
in~\eqref{eq:hhforD} we get the morphisms 
\eq\label{eq:HHEU3}
&&\C_\Delta\to \Omega_{X\times X}\lltens[\shd_{X\times X}](\shm\detens\RDD\shm)\to\omega_\Delta.
\eneq
For $\shm\in\Derb_\coh(\shd_X)$, we set
\eqn
&&\HK(\shm)\eqdot  \Omega_{X\times X}\lltens[\shd_{X\times X}](\shm\detens\RDD\shm).
\eneqn
Then $\HK(\shm)$ is a trace kernel by~\eqref{eq:HHEU3} and 
$\mueu_M(\HK(\shm))$ is supported by $\chv(\shm)$ by Theorem~\ref{th:ssinchar1}.

\begin{proposition}
The Hochschild class of $\shm$ is the microlocal Euler class of the trace kernel associated to $\shm$, that is,
$\hh_\she(\shm)=\mueu_X(\HK(\shm))$ in $H^0_{\chv(\shm)}(T^*X;\opb{\pi}\omega_X)$.
\end{proposition}
${}$

\section[Elliptic pairs]{Microlocal Euler class of elliptic pairs}

Let $X$ be a complex manifold,  $\shm$ an object of $\Derb_\coh(\shd_X)$ and $G$ an object 
of $\Derb_\Rc(\C_X)$. 
The pair $(\shm,G)$ is called an elliptic pair in~\cite{ScSn94} if $\chv(\shm)\cap\SSi(G)\subset T^*_XX$. From now on, we assume that  $(\shm,G)$ is an elliptic pair.
We set
\eq\label{eq:ellpairs2}
&&\HK(\shm,G)\eqdot
   \Omega_{X\times X}\ltens[\shd_{X\times X}]\bl(\shm\tens G)\detens(\RDD\shm\tens\RD_X'G)\br.
\eneq
It follows from the preceding results that $\HK(\shm,G)$ is a trace kernel and 
\eq
&&\mueu_X\bl\HK(\shm,G)\br=\mueu_X(\shm)\star\mueu_X(G).
\eneq
%We set
%\eq
%&&\GSol(\shm,F)\eqdot \RHom[\shd_X](\shm\tens F,\sho_X),\label{eq:SolMF}\\
%&&\DR(\shm,F)\eqdot \rsect(X;\Omega_X\lltens[\shd_{X}]\shm\tens F)\,[d_X].\label{eq:DRMF}
%\eneq

Applying Corollary~\ref{cor:petrowshv}~(a),  we get the natural isomorphism
\eq\label{eq:petrovski}
&&\rhom[\shd_X](\shm,\RD_X'G\tens\sho_X)\isoto\rhom[\shd_X](\shm\tens G,\sho_X).
\eneq
Assume moreover that $\Supp(\shm)\cap\Supp(G)$ is compact. Applying Corollary~\ref{cor:petrowshv}~(b), we get that
the cohomology of the complex 
\eqn
&&\GSol(\shm\tens G)\eqdot\RHom[\shd_X](\shm\tens G,\sho_X)
\eneqn 
is finite dimensional. Moreover
\eqn
&&\rsect(X\times X;\HK(\shm,G))\simeq \GSol(\shm\tens G)\tens\GSol(\shm\tens G)^*.
\eneqn
Applying Proposition~\ref{pro:convmueu}, we get
\eqn
\chi\bl\rhom[\shd_X](\shm\tens G,\sho_X)\br&=&\int_X(\hh_\she(\shm)\star\mueu_X(G))\vert_X\\
&=&\int_{T^*X}(\hh_\she(\shm)\cup\mueu_X(G)).
\eneqn
This formula has many applications, as far as one is able to calculate 
$\mueu_X(\shm)$. 

Assume that $\shm$ is endowed with a good filtration and $\chv(\shm)\subset\Lambda$. Set
\eqn
&&\tw \gr\shm\eqdot \sho_{T^*X}\tens[\opb{\pi}\gr\shd_X]\opb{\pi}\gr\shm\\
&&\sigma_\Lambda(\shm)=\Ch_\Lambda(\tw\gr\shm)\in\bigoplus_j H^{2j}_\Lambda(T^*X;\C_{T^*X}),\\
&&\mu\Ch_\Lambda(\shm)=\sigma_\Lambda(\shm)\cup \spb{\pi}\Td_X(T^*X)\mbox{ for a left $\shd$-module},\\
&&\mu\Ch_\Lambda(\shm)=\sigma_\Lambda(\shm)\cup \spb{\pi}\Td_X(TX)\mbox{ for a right $\shd$-module},
\eneqn
where $\Ch$ is the Chern character and $\Td$ is the Todd class.
Note that $\mu\Ch$ commutes with proper direct images (Laumon's version of the RR theorem for $\shd$-modules) and non characteristic inverse images. 
In~\cite{ScSn94} we made the conjecture that
\eqn
&&\mueu_\Lambda(\shm)=[\mu\Ch_\Lambda(\shm)]^{2 d_X}
\eneqn
This conjecture has been proved 
in~\cite{BNT02}  by Bressler-Nest-Tsygan and generalized  in~\cite{BGNT07}.

\begin{example}
(i) If $X$ is a complex compact manifold, one recovers the Riemann-Roch theorem: one takes $G=\C_X$ and if $\shf$ is a coherent $\sho_X$-module, one sets $\shm=\shd_X\tens[\sho_X]\shf$. 

\vspace{0.3ex}\noindent
(ii) If $M$ is a compact real analytic manifold and $X$ is a complexification of $M$, one recovers the Atiyah-Singer theorem by choosing $G=\RD_X'\C_M$. 
\end{example}

${}$

\chapter[Indsheaves]{Ind-sheaves and applications to $\shd$-modules}

%\begin{abstract}
{\bf Abstract.}
I will first recall the constructions of~\cite{KS96,KS01} of the sheaves of temperate or Whitney  holomorphic functions. These are not sheaves on the usual topology, but sheaves on the subanalytic site or better, ind-sheaves. Then I will explain how these objects appear naturally in the study of irregular holonomic $\shd$-modules.
%\end{abstract}

\section{Ind-sheaves}

\subsubsection*{Ind-objects}
References are made to~\cite{SGA4} or to~\cite{KS06} for an exposition. We keep the notations of the preceding lectures.

Let $\shc$ be an abelian category (in a given universe $\shu$). One denotes by $\shc^{\wedge,add}$ the big category of
additive  functors  from $\shc^\rop$ to $\md[\Z]$. 
This big category is abelian and the functor 
$h^\wedge\cl \shc \to \shc^\wedge$
makes $\shc$ a full abelian subcategory of $\shc^{\wedge,add}$. This
functor is left exact, but not exact in general.

An ind-object in $\shc$ is an object 
$A\in \shc^\wedge$ which is isomorphic
to $\inddlim \alpha$ for some functor $\alpha\colon I \to \shc$ with $I$
filtrant and small. One denotes by $\indc$ the full additive  subcategory of 
$\shc^{\wedge,add}$ consisting of ind-objects. 

\begin{theorem}
\bnum
\item 
The category $\indc$ is abelian.
\item 
The natural functor $\shc\to \indc$ is fully faithful and exact
and the natural functor $\indc\to \shc^{\wedge,add}$ is fully faithful
and left exact. 
\item 
The category $\indc$ admits exact small filtrant inductive limits
and the functor $\indc\to \shc^{\wedge,add}$ commutes with such limits.
\item 
Assume that $\shc$ admits small projective limits. 
Then the category $\indc$ admits small 
projective limits, and the functor $\shc\to \indc$ commutes with such limits.
\enum
\end{theorem}

\begin{example}
Assume that $\cor$ is a field and denote by $\mdf[\cor]$ the category of finite dimensional $\cor$-vector spaces.
Let $\II[\cor]$ denote the category of ind-objects of $\md[\cor]$. Define 
$\beta\colon\md[\cor]\to\II[\cor]$ by 
setting $\beta(V)=\sinddlim W$, where $W$ ranges over the
family of finite-dimensional vector subspaces of $V$. 
In other words, $\beta(V)$ is the functor from $\md[\cor]^\rop$ to $\md[\Z]$
given by $M\mapsto\indlim[W]\Hom[k](M,W)$. Therefore,
\eqn
\indlim[W\subset V,{W\in\mdf[\cor]}]\Hom[\cor](L,W)
    &\simeq& \Hom[{\II[\cor]}](L,\inddlim[W\subset V,{W\in\mdf[\cor]}]W)\\
    &=& \Hom[{\II[\cor]}](L,\beta( V)).
\eneqn
If $V$ is infinite-dimensional, $\beta(V)$ is not representable in $\md[\cor]$.  Moreover, 
$\Hom[{\II[\cor]}](\cor,V/\beta(V))\simeq 0$.
\end{example}
It is proved in~\cite{KS06} that the category $\indc$ for $\shc=\md[\cor]$ does not have enough injectives.

\begin{definition}
An object $A\in\indc$ is quasi-injective if the functor $\Hom[\indc](\scbul,A)$ is exact on the category $\shc$. 
\end{definition}
It is proved in loc.\ cit.\  that if $\shc$ has enough injectives, then $\indc$ has enough quasi-injectives.

\subsubsection*{Ind-sheaves}
References are made to~\cite{KS01}.

Let $X$ be a locally compact space countable at infinity.
Recall that $\md[\cor_X]$ denotes the abelian category of sheaves of $\cor$-modules on $X$.
We denote by $\mdcp[\cor_X]$ the full subcategory consisting of sheaves with compact support. 
We set for short:
\eqn
&&\II[\cor_X]:= {\rm Ind}(\mdcp[\cor_X])
\eneqn
and call an object of this category an indsheaf on $X$. 

\begin{theorem}
The prestack $U\mapsto\II[\cor_U]$, $U$ open in $X$ is a stack.
\end{theorem}

The following example explains why we have considered sheaves with compact supports. 
\begin{example}
Let $X=\R$, let $F=\cor_X$, $G_n=\cor_{[n,+\infty[}$, $G=\inddlim[n]G_n$.
Then $G\vert_U=0$ in ${\rm Ind}(\md[\cor_U])$ 
for any relatively compact open subset $U$ of $X$.
On the other hand,
$\Hom[{{\rm Ind}(\md[\cor_X])}](\cor_X,G)\simeq\indlim[n]\Hom[\cor_X](\cor_X,G_n)\simeq \cor$.
\end{example}

We have two pairs $(\alpha_X,\iota_X)$ and $(\beta_X,\alpha_X)$  of adjoint functors 
\eqn
&& \xymatrix{
{\md[\cor_X]}\ar@<-1.5ex>[rr]_-{\beta_X}\ar@<1.5ex>[rr]^-{\iota_X}&&{\II[\cor_X]}\ar[ll]|-{\alpha_X}.
}
\eneqn
The functor $\iota_X$ is the natural one. If $F$ has compact support, $\iota_X(F)=F$ after identifying a category $\shc$ to a full subcategory of $\indc$. The functor $\alpha_X$ associates 
$\indlim[i] F_i$ ($F_i\in\mdcp[\cor_X]$, $i\in I$, $I$ small and filtrant) to the object $\inddlim[i] F_i$. If $\cor$ is a field, $\beta_X(F)$ is the functor $G\mapsto\sect(X;H^0(\RD_X'G)\tens F)$. 
\begin{itemize}
\item
$\iota_X$ is exact, fully faithful, and commutes with $\sprolim$, 
\item $\alpha_X$ is exact and commutes with $\sprolim$ and $\sindlim$,
\item $\beta_X$ is right exact, fully faithful and commutes with $\sindlim$,
\item $\alpha_X$ is left  adjoint to $\iota_X$,
\item $\alpha_X$ is right adjoint to $\beta_X$,
\item $\alpha_X\circ\iota_X\simeq \id_{\md[\cor_X]}$ 
and $\alpha_X\circ\beta_X\simeq \id_{\md[\cor_X]}$.
\end{itemize}
\begin{example}
Let $U\subset X$ be an open subset, $S\subset X$ a closed subset. Then
\eqn
&&\beta_X(\cor_U)\simeq\inddlim[V]\cor_V,\, V\mbox{ open }, V\subset \subset U,\\
&&\beta_X(\cor_S)\simeq\inddlim[V]\cor_{\ol V},\, V\mbox{ open }, S\subset V.
\eneqn

Let $a\in X$ and consider the skyscraper sheaf $\cor_{\{a\}}$. Then $\beta_X(\cor_{\{a\}})\to\cor_{\{a\}}$ is an epimorphism in $\II[\cor_X]$ and defining $N_a$ by the exact sequence:
\eqn
0\to N_a\to \beta_X(\cor_{\{a\}})\to\cor_{\{a\}}\to0
\eneqn
we get that $\Hom[{\II[\cor_X]}](\cor_U,N_a)\simeq 0$ for all open neighborhood $U$ of $a$. 
\end{example}
We shall not recall here the construction of the derived category of indsheaves,  nor the six operations on such ``sheaves''.

\section{Sheaves on the subanalytic site}
The subanalytic site was introduced in~\cite[Chapt~7]{KS01} and the results on sheaves on this site were obtained as particular cases of more general results on indsheaves, which makes the reading not so easy. A direct and more elementary study of
sheaves on the subanalytic site is performed in~\cite{Pr08,Pr12}.

Let $M$ be a real analytic manifold. 
One denotes by $\rc[\cor_M]$ the abelian category of $\R$-constructible sheaves on $M$
and by $\rcc[\cor_M]$ the full subcategory consisting of sheaves with compact support. There is an equivalence $\Derb(\rc[\cor_M])\simeq\Derb_{\Rc}(\cor_{M})$ where this last category 
is the full triangulated subcategory of 
$\Derb(\cor_{M})$ consisting of  $\R$-constructible sheaves. (This classical result has first been proved by Kashiwara~\cite{Ka84}.)

We denote by $\Op_M$ the category whose objects are the 
open subsets of $M$ and the morphisms are the inclusions of open
subsets. One defines a Grothendieck topology on $\Op_M$ by 
deciding that a family  $\{U_i\}_{i\in I}$
of subobjects of $U\in\Op_M$ is a covering of $U$ if it is a covering
in the usual sense.
\begin{definition}
Denote by $\Op_{\Msa}$ the full subcategory of $\Op_M$ consisting of 
 subanalytic and relatively compact open subsets. 
 The site $\Msa$ is obtained by  deciding that a family  $\{U_i\}_{i\in I}$
of subobjects of $U\in\Op_{\Msa}$ is a covering of $U$ 
if there exists a finite subset $J\subset I$ such that $\bigcup_{j\in J}U_j=U$. 
\end{definition}
Let us denote by 
\eq
&&\rhosa\cl M\to\Msa
\eneq
the natural morphism of sites. Here again, we have two pairs of adjoint
functors $(\opb{\rhosa},\oim{\rhosa})$ and $(\eim{\rhosa},\opb{\rhosa})$ :

\eqn
&& \xymatrix{
{\md[\cor_M]}\ar@<-1.5ex>[rr]_-{\eim{\rhosa}}\ar@<1.5ex>[rr]^-{\oim{\rhosa}}&&{\md[\cor_\Msa]}\ar[ll]|-{\opb{\rhosa}}.
}
\eneqn
For $F\in \md[\cor_M]$, $\eim{\rhosa}F$ is the sheaf associated to the
presheaf $U\mapsto F(\overline U)$, $U\in \Op_{\Msa}$. 

\begin{proposition}
The restriction of the functor $\oim{\rhosa}$ to the category $\rc[\cor_M]$
is exact and fully faithful.
\end{proposition}
By this result, we shall consider the category $\rc[\cor_M]$ as  a full subcategory of 
$\md[\cor_M]$ as well as a full subcategory of $\md[\cor_{\Msa}]$.
Set 
$$\IIrc[\cor_M]={\rm Ind}(\rcc[\cor_M]).$$

\begin{theorem}
The natural functor $\alpha_{\Msa}\cl\IIrc[\cor_M]\to \md[\cor_{\Msa}]$ is an equivalence of categories.
\end{theorem}
In other words, ind-$\R$-constructible sheaves are ``usual sheaves''
on the subanalytic site.
By this result, the embedding $\rcc[\cor_M]\into\mdcp[\cor_M]$ gives a functor
$I_M\cl \md[\cor_{\Msa}]\to\II[\cor_M]$.
Hence, we have a quasi-commutative diagram of categories
\eq\label{dia00}
\xymatrix{
\md[\cor_M]\ar[r]^-{\iota_M}&\II[\cor_M]\\
\mdrc[\cor_M]\ar[u]\ar[r]^-{\oim{\rhosa}}&\md[\cor_{\Msa}]\ar[u]_-{I_M}
}
\eneq
in which all arrows are exact and fully faithful.
One shall be aware that  the diagram:
\eq\label{dia1}
\xymatrix{
\md[\cor_M]\ar[r]^-{\iota_M}\ar[rd]_-{\oim{\rhosa}}&\II[\cor_M]\\
\ar@{}[ru]|(.66){NC}&\md[\cor_{\Msa}]\ar[u]_-{I_M}
}
\eneq
is not commutative. Moreover, $\iota_M$ is exact and $\oim{\rhosa}$ is not right exact in general. 

One denotes by $\sinddlim$ the inductive limit in the category $\md[\cor_{\Msa}]$. 
One shall be aware that the functor $I_M$   commutes with inductive limits but $\oim{\rhosa}$ does not.

\section{Moderate and formal cohomology}

From now on, $\cor=\C$. As usual, we denote by 
$\shc_M^\infty$ (resp.\ $\shc_M^\omega$) the sheaf of complex functions
of class $C^\infty$ (resp.\ real analytic), by $\shd b_M$ 
(resp.\ $\shb_M$) 
the sheaf of Schwartz's distributions (resp.\ Sato's hyperfunctions),
and by $\shd_M$ the sheaf of analytic finite-order differential operators.
We also use the notation $\sha_M=\shc_M^\omega.$

\begin{definition} Let $U\in\Op_\Msa$ and 
let $f\in \shc^{\infty}_M(U)$. One says that $f$ has
{\it  polynomial growth} at $p\in M$
if it satisfies the following condition.
For a local coordinate system
$(x_1,\dots,x_n)$ around $p$, there exist
a sufficiently small compact neighborhood $K$ of $p$
and a positive integer $N$
such that
\eq
&\sup_{x\in K\cap U}\big(\dist(x,K\setminus U)\big)^N\vert f(x)\vert
<\infty\,.&
\eneq
It is obvious that $f$ has polynomial growth at any point of $U$.
We say that $f$ is  temperate at $p$ if all its derivatives
have polynomial growth at $p$. We say that $f$ is temperate 
if it is temperate at any point.
\end{definition}

For $U\in\Op_\Msa$,
denote by $\Cinft_M(U)$ the subspace
of $\shc^{\infty}_M(U)$ consisting of tempered functions. 

Denote by $\Dbt_M(U)$ the space of tempered distributions on $U$, 
defined by the exact sequence
\eqn
&&0\to\sect_{M\setminus U}(M;\Db_M)\to\sect(M;\shd b_M)\to\Dbt_M(U)\to 0.
\eneqn                     

Using  Lojasiewicz's inequalities \cite{Lo61} (see also~\cite{Ma67}), one easily proves
that 
\begin{itemize}
\item
the presheaf  $U\mapsto \Cinft_M(U)$ is a sheaf on  $\Msa$,
\item
the presheaf  $U\mapsto \Dbt_M(U)$ is a sheaf on $\Msa$. 
\end{itemize}
One denotes by $\Cinft_{\Msa}$ the first one and calls it the sheaf of 
 temperate $C^\infty$-functions. 
One denotes by $\Dbt_{\Msa}$ the second one and calls it the sheaf of 
 temperate distributions. 
Let $F\in\Derb_\Rc(\C_M)$. One has the isomorphism
\eq\label{eq:thom}
&&\opb{\rhosa}\rhom(F,\Dbt_{\Msa})\simeq \thom(F,\Db_M)
\eneq
where the right-hand side was defined by Kashiwara  as the main tool for his proof of the Riemann-Hilbert correspondence in~\cite{Ka80,Ka84}.

For a closed subanalytic subset $S$ in $M$, 
denote by $\shi^\infty_{M,S}$ the subsheaf 
of $\shc^{\infty}_M$ consisting of functions which vanish up to 
infinite order on $S$.
In~\cite{KS96}, one introduces the sheaf:
\eqn
\C_U\wtens \shc^{\infty}_M&:=&V\mapsto\sect(V;\shi^\infty_{V,V\setminus U})
\eneqn
and shows how to extend this construction and define an exact functor 
$\scbul\wtens \shc^{\infty}_M$ on $\mdrc[\C_M]$.
One denotes by $\shc_M^{\infty,\rm w}$ the sheaf on $\Msa$ given by 
\eqn
&&\shc_M^{\infty,\rm w}(U)=\sect(M;H^0(\RD'_M\cor_U)\wtens \shc_M^{\infty}), U\in \Op_{\Msa}.
\eneqn
If $\RD'_M\C_U\simeq\C_{\ol U}$, $\shc_M^{\infty,\rm w}(U)$ is the space of Whitney functions on $U$, that is the quotient $\shc^\infty(M)/\shi^\infty_{M,M\setminus U}$. 
It is thus natural to call $\shc_M^{\infty,\rm w}$the sheaf of Whitney $C^\infty$-functions on $\Msa$. 

 Note that the sheaf $\oim{\rhosa}\shd_{M}$ does not operate on
the sheaves  $\Cinft_M$, $\shd b^t_M$, $\C_M^{\infty,\rm w}$ but  $\eim{\rhosa}\shd_{M}$ does.

\vspace{3ex}
Now let $X$ be a {\em complex} manifold. We still denote by $X$ 
the real underlying manifold and we denote by $\overline X$ the
complex manifold conjugate to $X$. 
One defines the sheaf of temperate holomorphic functions 
$\Ot[\Xsa]$ as the Dolbeault complex with coefficients in  
$\Cinft_{\Xsa}$. More precisely
\eq\label{eq:defOt1}
\Ot[\Xsa]&=&\rhom[\eim{\rhosa}\shd_{\overline X}](\eim{\rhosa}\sho_{\overline X},\Cinft_{\Xsa}).
\eneq
One proves the isomorphism
\eq\label{eq:defOt2}
\Ot[\Xsa]&\simeq&\rhom[\eim{\rhosa}\shd_{\overline X}](\eim{\rhosa}\sho_{\overline X},\Dbt_{\Xsa}).
\eneq
Similarly, one defines the sheaf 
\eq
&&\Ow[\Xsa]=\rhom[\eim{\rhosa}\shd_{\overline X}](\eim{\rhosa}\sho_{\overline X},\shc^{\infty,\rm w}_{\Xsa}).
\eneq
Note that the objects $\Ot[\Xsa]$ and   $\Ow[\Xsa]$ are  not
concentrated in degree zero in dimension $>1$. Indeed, with the subanalytic topology, only finite coverings 
are allowed. If one considers for example the open set $U\subset\C^n$, the difference of a open ball of radius $R>0$ and a closed ball of radius $r$ with  $0<r<R$, then the Dolbeault complex will not be exact after
 any finite covering. 
 For the same reason,  the sheaf $\roim{\rhosa}\sho_X$ is not concentrated in degree zero in dimension $>1$. 
 
Therefore, we shall better consider indsheaves and we shall embed the category 
$\Derb(\C_{\Xsa})$ into the category $\Derb(\II[\C_X])$ by the exact functor $I_X$. Hence we consider subanalytic sheaves as indsheaves. In the category  $\Derb(\II[\C_X])$ we have thus the morphisms of sheaves
\eqn
&&\Oww\to\Ow\to\Ot\to\sho_X.
\eneqn
Here $\Ow$ and $\Ot$ are the images of $\Ow[\Xsa]$ and  $\Ot[\Xsa]$ by the functor $I_X$ (there are still not concentrated in degree $0$), we have kept the same notation for $\sho_X$ and its image in $\md[{\II[\C_X]}]$ by the functor $\iota_X$,  and we have set
\eqn
&&\Oww\eqdot\beta_X(\sho_X).
\eneqn
We call $\Ow$ and $\Ot$ the sheaves of temperate and Whitney holomorphic functions, respectively. 
  
\begin{example}
Let $Z$ be a closed complex analytic subset of the complex manifold $X$. We have the isomorphisms
\eqn
&&\alpha_X\rhom[{\II[\C_X]}](\RD'\C_Z,\Oww)\simeq \sho_X\vert_Z,\\
&&\alpha_X\rhom[{\II[\C_X]}](\RD'\C_Z,\Ow)\simeq
 \sho_X{\widehat \vert}_Z\mbox{ (formal completion along $Z$),}\\
&&\alpha_X\rhom[{\II[\C_X]}](\C_Z,\Ot)\simeq\rsect_{[Z]}(\sho_X)\mbox{ (algebraic cohomology)},\\
&&\alpha_X\rhom[{\II[\C_X]}](\C_Z,\sho_X)\simeq\rsect_{Z}(\sho_X).
\eneqn
\end{example}
\begin{example}
let $M$ be a real analytic manifold and $X$ a complexification of $M$. We have the isomorphisms
\eqn
&&\alpha_X\rhom[{\II[\C_X]}](\RD'\C_M,\Oww)\simeq \sha_M,\\
&&\alpha_X\rhom[{\II[\C_X]}](\RD'\C_M,\Ow)\simeq \shc_M^\infty,\\
&&\alpha_X\rhom[{\II[\C_X]}](\RD'\C_M,\Ot)\simeq\Db_M,\\
&&\alpha_X\rhom[{\II[\C_X]}](\RD'\C_M,\sho_X)\simeq\shb_M.
\eneqn
\end{example}
Notice that with this approach, the sheaf $\Db_M$ of Schwartz's distributions is constructed similarly as the sheaf of Sato's hyperfunctions. In particular,  functional analysis is not used in the construction.

\begin{remark}
The subanalytic topology allows us to define functions whose growth at the boundary is bounded by some power of the inverse of the distance to the boundary, but not to make precise this power. In order to define such sheaves, we have recently defined with S.~Guillermou in~\cite{GS13} the linear subanalytic topology $\Msal$ on a real analytic manifold $M$. The open sets of this topology are those of $\Msa$, namely $\Op_\Msa$, but there are less coverings. Roughly speaking, a finite covering $\{U_i\}_{i\in I}$ 
is a linear covering of $U=\bigcup_iU_i$ if there is a constant $C$ such that for any $x\in M$
\eq\label{eq:resit}
&&d(x,M\setminus\bigcup_{i\in I}U_i)\leq C\cdot
\max_{i\in I}  d(x,M \setminus U_i).
\eneq
Here $d$ is a distance on $M$ which is locally equivalent to the Euclidian distance on $\R^n$.
One proves that the family of linear coverings satisfies the axioms
of Groth\-endieck topologies. One denotes by $\Msal$ the site so defined and by 
$\rhosal\cl \Msa\to \Msal$
the natural morphism of sites. One of the main results of the theory is that 
functor $\roim{\rhosal}\cl\RD^+(\cor_\Msa)\to\RD^+(\cor_\Msal)$ admits a right adjoint
 $\epb{\rhosal}\cl\RD^+(\cor_\Msal)\to\RD^+(\cor_\Msa)$. Moreover, if $U\in\Op_\Msa$ has Lipschitz boundary, then $\roim{\rhosal}\C_U$ is concentrated in degree $0$. It follows that if $F$ is a presheaf 
 on $\Msa$ such that the sequence $0\to F(U_1\cup U_2)\to F(U_1)\oplus F(U_2)\to F(U_1\cap U_2)\to 0$ is exact for any linear covering $(U_1,U_2)$ of $U_1\cup U_2$, then there exists  
 $F\in\RD^+(\cor_\Msa)$ such that $\rsect(U;F)\simeq F(U)$ for all $U\in\Op_\Msa$ with Lipchitz boundaries.

 This topology allows us to define the subsheaf $\Cin[s]_\Msal$ 
of $\Cinf_\Msal$ consisting of  functions tempered of order $s$. On a complex manifold $X$ we may thus endow the sheaf $\Ot[\Xsa]$ with a natural filtration (in the derived sense). We refer to loc.\ cit.\ more more details.
\end{remark}

\section{Applications to $\shd$-modules I}

Let us show on an example extracted of~\cite{KS03} the possible role of the sheaf $\sho_X^t$ 
in the study of irregular holonomic $\shd$-modules. 

Let $X$ be a complex manifold and let $\shm$ be a holonomic  $\shd$-module. 
We set for short
\eqn
\Solo(\shm)&=& \hom[\shd_X](\shm,\sho_X),\\
\Solot(\shm)&=&\hom[\beta_X\shd_X](\beta_X\shm,\sho_X^t).
\eneqn
We shall compare these two objects in a simple example in which $\shm$  is not regular.
Let $X=\C$ endowed with the holomorphic coordinate $z$ and let $P=z^2\partial_z+1$. 
We consider the $\shd_X$-module $\shm\eqdot\shd_X\exp(1/z)\simeq\shd_X/\shd_X\cdot P$.

Notice first that $\sho_X^t$ is concentrated in degree $0$ (since
$\dim X=1$) and it is a sub-indsheaf of $\sho_X$. It follows that the
morphism $\Solot(\shm)\to\Solo(\shm)$ is a monomorphism.
Moreover,
\eqn
&&\Solo(\shm)\simeq\C_{X,X\setminus \{0\}}\cdot\exp(1/z).
\eneqn
It follows that for $V\subset X$ a connected open subset, we have
$\sect(V;\Solot(\shm))\neq 0$ if and only if $V\subset X\setminus \{0\}$ and
$\exp(1/z)\vert_V$ is tempered.

Let $\bar B_\epsilon$ denote the closed ball with center $(\epsilon,0)$
and radius $\epsilon$ and set
$U_\epsilon=X\setminus \bar B_\epsilon$. 

Then one proves  that
$\exp(1/z)$ is temperate (in a neighborhood of $0$)
 on an open subanalytic subset $V\subset X\setminus\{0\}$ if and only
if ${\rm Re}(1/z)$ is bounded on $V$, that is, if and only if
$V\subset U_\epsilon$ for some $\epsilon>0$. We get

\begin{proposition}\label{th:exholirr}
One has the isomorphism
\eq\label{eq:solt=indlim}
&&
\inddlim[\epsilon>0]\C_{XU_\epsilon}\isoto\Solot(\shm).
\eneq
\end{proposition}
Unfortunately, the functor $\Solt$ (as well as its derived functor) is not fully faithful since 
the $\shd$-modules  $\shm\eqdot\shd_X\exp(1/z)$ and $\shn\eqdot\shd_X\exp(2/z)$ have the same 
indsheaves of temperate holomorphic solutions although they are not isomorphic.

Proposition~\ref{th:exholirr} has been generalized to the study of holonomic modules in 
dimension one in~\cite{Mr09}.

\section{Applications to $\shd$-modules II}
For $F\in\Derb_{\Rc}(\C_X)$, set (see~\eqref{eq:thom}):
\eqn
F\wtens\sho_X&\eqdot& \rhom[\shd_{\overline X}](\sho_{\overline X},F\wtens\shc^\infty_X),\\
\thom(F,\sho_X)&\eqdot &
\rhom[\shd_{\overline X}](\sho_{\overline X},\thom(F,\Db_X)).
\eneqn

Let $F\in\Derb_{\Rc}(\C_X)$ and $\shm\in\Derb_{\coh}(\shd)$. Recall that we have set
$\omega_X^\hol\eqdot \Ov[X][d_X]$. 
Set for short
\eqn
&&W(\shm,F)\eqdot \rhom[\shd](\shm,F\wtens\sho_X), \\
&&T(F,\shm) \eqdot\thom(F,\omega_X^\hol)\ltens[\shd]\shm.
\eneqn
There is a natural morphism
\begin{equation}\label{eq:pre-paring}
W(\shm,F) \tens T(F,\shm) \to \omega_X^\hol,
\end{equation} 
functorial in $F$ and $\shm$.
For $G\in\Derb_{\Rc}(\C_X)$ one gets a pairing
\begin{align}
\label{eq:pairing}
\RHom(G,W(\shm,F)) &\tens \rsect_c(X;G\tens T(F,\shm)) \\
\notag
&\to \rsect_c(X;W(\shm,F) \tens T(F,\shm)) \\
\notag
&\to \rsect_c(X;\omega_X^\hol)\to \C.
\end{align}
Denote by $\Derb(FN)$ the derived category of the quasi-abelian category of  Fr\'echet nuclear $\C$-vector spaces and define similarly the category $\Derb(DFN)$, where now DFN stands for ``dual of Fr\'echet nuclear''.

\begin{theorem}\label{thm:KSduality} {\rm (\cite[Theorem 6.1]{KS96})}\,
Let $F,G\in\Derb_{\Rc}(\C_X)$ and $\shm\in\Derb_{\coh}(\shd)$. Then the two
complexes
\eqn
&&\RHom(G,W(\shm,F))\in \Derb(FN)\mbox{ and }
\rsect_c(X;G\tens T(F,\shm)) \in \Derb(DFN)
\eneqn
are dual to each other through \eqref{eq:pairing}, functorially in $F,G$ and $\shm$.
\end{theorem}
Now we assume that $\shm\in\Derb_\hol(\shd_X)$ and we consider the following assertions.
\begin{itemize}
\item[(a)]
$W(\shm,F)=\rhom[\shd](\shm,F\wtens\sho_X)$ is $\R$-constructible,
\item[(b)]
$T(F,\shm)=\thom(F,\omega_X^\hol)\ltens[\shd]\shm$ is $\R$-constructible,
\item[(c)]
the two complexes in (a) and (b) are dual to each other in the category $\Derb_\Rc(\C_X)$, that is,
$W(\shm,F) \simeq \RD_X T(F,\shm)$.
\end{itemize}
It was conjectured in~\cite{KS03} that (b) is always satisfied.  
Based on the work of Mochizuki~\cite{Mo09} (see also \cite{Ke10,Ke11,Sb12}),  partial 
results in this direction have been obtained in~\cite{Mr13}.

On the other hand, one deduces easily  from 
Theorem~\ref{thm:KSduality} that (a) and (b) are equivalent and imply (c). 
Finally, it follows immediately from~\cite{Ka78,Ka84} that (b), hence (a) and (c), are true 
when  $F\in\Derb_{\Cc}(\C_X)$.

\begin{corollary}\label{cor:KSduality} 
Assume that $F\in\Derb_{\Cc}(\C_X)$ and $X$ is compact.
Then the complexes $\rsect(X;W(\shm,F))$ and $\rsect(X;T(F,\shm))$ have finite-dimensional cohomology
 and~\eqref{eq:pairing} induces a perfect pairing for all $i\in\Z$
\[
H^{-i}\rsect(X;W(\shm,F)) \tens H^i\rsect(X;T(F,\shm))\to\C,
\]
functorial in $F$ and $\shm$.
\end{corollary}
In~\cite{BE04}, S.~Bloch and H.~Esnault prove directly a similar result on an algebraic curve $X$ when assuming that $\shm$ is a meromorphic connection with poles on a divisor $D$. They interpret the duality pairing by considering sections of the type 
$\gamma\tens\epsilon$, where $\gamma$ is a  cycle  with boundary on $D$ and $\epsilon$ is a horizontal section of the connection 
on $\gamma$ with exponential decay on $D$. Their work has been extended to  higher dimension  by M.~Hien~\cite{H09}.

It would be interesting to make a link with these results and Corollary~\ref{cor:KSduality}.

\providecommand{\bysame}{\leavevmode\hbox to3em{\hrulefill}\thinspace}

\parbox[t]{21em}
{\scriptsize{
\noindent

Pierre Schapira\\
Institut de Math{\'e}matiques,
Universit{\'e} Pierre et Marie Curie\\
and Mathematics Research Unit, University of Luxemburg\\
%4 Place Jussieu, 7505 Paris\\
e-mail: schapira@math.jussieu.fr\\
http://www.math.jussieu.fr/\textasciitilde schapira/
}}

\end{document}